\documentclass[preprint,12pt]{elsarticle}

\usepackage{amssymb}
\usepackage{amsmath}

\journal{Journal of Computational Physics}

\begin{document}

\begin{frontmatter}



\title{Stochastic estimation of Green's functions with application to diffusion and advection-diffusion-reaction problems}


\author[inst1]{Russell G. Keanini}

\affiliation[inst1]{organization={Department of Mechanical Engineering and Engineering Science},
            addressline={9201 University City Blvd.}, 
            city={Charlotte},
            postcode={28223}, 
            state={North Carolina},
            country={United States}}

\author[inst2]{Jerry Dahlberg}
\affiliation[inst2]{organization={University of Tennessee Space Institute}, city={Huntsville},
postcode={35649},
state={Alabama},
country={United States}}

\author[inst1]{Philip Brown}
\author[inst1]{Mehdi Morovati}
\author[inst1]{Hamidreza Moradi}
\author[inst3]{Donald Jacobs}
\affiliation[inst3]{organization={Department of Physics and Optical Science},
            addressline={9201 University City Blvd.}, 
            city={Charlotte},
            postcode={28223}, 
            state={North Carolina},
            country={United States}}
\author[inst1]{Peter T. Tkacik}

\begin{abstract}
A stochastic method is described for estimating Green's functions (GF's), appropriate to linear advection-diffusion-reaction transport problems, evolving in arbitrary geometries. By allowing straightforward construction of approximate, though high-accuracy GF's, within any geometry, the technique solves the central challenge in obtaining Green's function solutions. In contrast to Monte Carlo solutions of individual transport problems, subject to specific sets of conditions and forcing, the proposed technique produces approximate GF's that can be used: a) to obtain (infinite) sets of solutions, subject to any combination of (random and deterministic) boundary, initial, and internal forcing, b) as high fidelity direct models in inverse problems, and c) as high quality process models in thermal and mass transport design, optimization, and process control problems.  The technique exploits an equivalence between the adjoint problem governing the transport problem Green's function, $ G \left( \mathbf{x}, t | \mathbf{x'} ,t' \right) , $ and the backward Kolmogorov problem governing the transition density, $ p \left( \mathbf{x}, t | \mathbf{x'} ,t' \right) ,$ of the stochastic process used in Green's function construction.  We address nonspecialists and report four contributions. First, a recipe is outlined for diagnosing when  stochastic Green's function estimation can be used, and for subsequently estimating the transition density and associated Green's function.  Second, a naive estimator for the transition density is proposed and tested. 
Third, Green's function estimation error produced by random walker absorption at Dirichlet boundaries is suppressed using a simple random walker splitting technique. 
Last, spatial discontinuity in estimated GF's, produced by the naive estimator, is suppressed using a simple area averaging method. The paper provides guidance on choosing key numerical parameters, and the technique is tested against two simple unsteady, linear heat conduction problems, and an unsteady groundwater dispersion problem, each having known, exact GF's.
\end{abstract}




\end{frontmatter}

\section{Introduction}
Considering a physical system whose dynamics are determined by a set of conservation laws and physical principles, finding a Green's function that reliably models the system's spatial and temporal response to forcing opens the door to a number of powerful capabilities. Focusing on thermal and chemical transport in complicated, irregular regions, 
the temporal thermal or chemical response at discrete locations, to an infinite array of space and time-dependent forcing, for example, can be readily determined \cite{cole, barton1989elements}.  Similarly, Green's functions can be incorporated in inverse problems, allowing, for example, extraction of thermal and  pollutant source locations and strengths, based on limited temperature and chemical measurements; see references and examples in \cite{inversegf, inversegf2, inversegf3}. Green's functions can also provide a basis for efficient thermal and mass transport design, optimization, and process control  \cite{procpaper1, procpaper2, procpaper3}.  In these applications, a reliable Green's function can replace approximate or empirical system models, improving the reliability and quality of a design or control process.    

Unfortunately, building a Green's function rests on three typically challenging tasks:

\vspace{0.3cm}

\noindent A) derive (or, when possible, look up) the adjoint equation governing the Green's function (GF);

\vspace{0.4cm}

\noindent B) simultaneously derive the integral solution for the evolving field variable, $ \eta \left( \mathbf{x}, t \right) $ - which we'll call the 'magic rule' \cite{barton1989elements} -  stated in terms of convolutions of the Green's function with initial and boundary conditions, as well as in-domain forcing terms, $ f(\mathbf{x}, t) ; $ and

\vspace{0.4cm}

\noindent C) find an appropriate Green's function, subject to the governing adjoint problem derived in steps A) and B). 

Most texts treating Green's function techniques do not illustrate the trial and error nature of the first two steps, leaving non-specialists with limited insight. We address this deficiency in Appendix C, simultaneously deriving, via trial and error, the adjoint equation and magic rule associated with a linear advection-diffusion problem.  

Importantly, most of the difficulty surrounding Green's function solutions is tied to task C). A number of texts illustrate this step, invariably focusing on solution of (given, non-derived) adjoint equations, stated in simple geometries, leading to analytical Green's function solutions \cite{cole,barton1989elements,economou}. For linear transport  problems in complex geometries, however, this step clearly requires a numerical attack. Unfortunately, little work has been reported on numerical construction of Green's function's, specifically, numerical solution of adjoint equations associated with transport problem-specific partial differential equations, such as Eqs. (\ref{mainequation1}) and (\ref{mainequation}) below. While boundary element (BEM) \cite{numericalconstructiongf,bemoverview,bembook}, Green element (GEM) \cite{gem1, gem2, gem3}, and finite element methods (FEM) \cite{fem1, fem2, fem3, fem4} are widely characterized as Green's function techniques, none of these explicitly solve problem-specific adjoint equations. As a consequence, these only allow whole-domain solutions, in contrast to the pointwise solutions obtained by the technique proposed here. 

We focus on construction of 
approximate GF's, appropriate to solution of advection-diffusion-reaction problems of the general form:
\begin{equation}\label{mainequation1}
\mathbf{D} \mathbf{:} \nabla^T \otimes \nabla \eta - \mathbf{v}\cdot\nabla\eta - \gamma\eta - \frac{\partial\eta}{\partial t} = -f(\mathbf{x}, t),
\end{equation}
subject to specified boundary and initial conditions. Here, $ \eta = \eta \left( \mathbf{x}, t \right) , $ is a scalar field property (often a density of some kind), subject to:  i) transport by diffusion/dispersion $ ( \mathbf{D} $ real-valued ) and/or ii) wave propagation $ ( \mathbf{D} $ complex-valued), iii) advective transport by a fluid velocity field, $ \mathbf{v} =\mathbf{v} \left( \mathbf{x}, t \right ), $ iv) depletion by first order chemical reaction, $ - \gamma \eta , $ and v) augmentation $ \left( f \left(\mathbf{x}, t \right) > 0 \right) $ or depletion $ \left( f \left(\mathbf{x}, t \right) < 0 \right) $ by a source/sink field, $ f \left(\mathbf{x}, t \right) . $ In addition, 
$ \mathbf{D} $ represents an anisotropic, often position- and time-dependent diffusion or dispersion tensor, $ \gamma $ is a solute-dependent chemical or radiological decay constant, where $ \mathbf{D} ,$ $ \mathbf{v} , $ and $ \gamma $ have no dependence on $ \eta , $  while $ f $ is at most, linearly dependent on $ \eta .$ Importantly, the dispersion/diffusivity, velocity and source/sink fields are known quantities. 

In situations where the dispersion/diffusion tensor, $ \mathbf{D} ,$ is diagonal - capturing diffusive or dispersive transport that's dominant in one, two, or three mutually orthogonal directions, the substitutions $ \eta \rightarrow \exp \left( -\gamma t \right) \eta  $ and  $ f \left( \mathbf{x}, t \right) \rightarrow \exp \left( -\gamma t \right) f \left( \mathbf{x}, t \right) $, can be used to transform Eq. (\ref{mainequation1}) into an equation of the form, 
\begin{equation}\label{mainequation}
\mathbf{D} \mathbf{:} {\nabla^T} \otimes \nabla  \eta - \mathbf{v}\cdot\nabla\eta - \frac{\partial\eta}{\partial t} = -f(\mathbf{x}, t),
\end{equation}
As discussed in the third validation test below, equations of the form (\ref{mainequation1}) arise, for example, in models of solute transport in aquifers \cite{sanskrityayn2018analytical}. Equation (\ref{mainequation}) appears in many areas of physics, biology, and engineering, including fluid mechanics (linearized Navier-Stokes equations) \cite{keanini2011green,jovanovic2001modeling}, drift-diffusion \cite{levich1962physicochemical} and reaction-diffusion problems \cite{williams1965combustion,grindrod1996theory}), quantum mechanics (Schrödinger equation) \cite{messiah2014quantum}, statistical mechanics (Fokker-Planck equation) \cite{landau1981course} and plasma physics (Vlasov equation) \cite{colonna2016boltzmann}. Wave equations \cite{whitham2011linear,drazin1989solitons} and Chapman-Kolmogorov equations \cite{gardiner2004handbook} can also be expressed in this form. 

Since 'backward time' stochastic construction of GF's provides the single-point response at a chosen 'response point', $ \left( \mathbf{x},t \right) ,$ as produced by any given set of boundary and initial conditions, we focus most of the discussion on backward time Green's function construction.  For testing and validation purposes, however, the last validation case uses forward time stochastic construction. As discussed in test 3, forward time construction: i) can be used when the diffusion/dispersion tensor is diagnonal: $ \mathbf{D} = \mathbf{\tilde{D}} \cdot \mathbf{I} , $ where $ \mathbf{\tilde{D}} $ is a tensor function of position and time, and $ \mathbf{I} $ is the identity matrix, and ii) provides the response at $ \left( \mathbf{x},t \right) ,$ as produced by an impulse at a specific 'impulse point', $ \left( \mathbf{x'},t' \right) .$ 

Focusing on backward time construction, the proposed approach launches a swarm of RW's backward in time from a chosen response point, $ \left( \mathbf{x},t \right) ,$ calculating the random, backward trajectories of each random walker, using a transport-problem-specific stochastic differential equation. The transition (probability) density, $ p\left( \mathbf{x} , t | \mathbf{x'} , t' \right),$ for observing the random walker at a chosen 'impulse point', $ \left( \mathbf{x}^{\prime} , t^{\prime} \right) , $ given that it started at $ \left( \mathbf{x},t \right) , $ is then estimated.  Under conditions where the transport problem's adjoint problem - governing evolution of the propagator, $ K\left( \mathbf{x} , t | \mathbf{x'} , t' \right),$ - has the same structure as the so-called backward Kolmogorov problem - governing evolution of $ p\left( \mathbf{x} , t | \mathbf{x'} , t' \right)$ - the estimated transition density  allows estimation of the propagator, and in turn, the Green's function.  

In overview, we first outline a step-by-step recipe for stochastically constructing GF's, appropriate to transport problems governed by Eqs. (\ref{mainequation1}) and (\ref{mainequation}). The proposed stochastic Green's function estimation technique is then contrasted against Monte Carlo approaches for (approximately) solving initial-boundary value problems (IBVP's) \cite{muller, sabelfeld, booth, sabelfeldsimonov, haji, mascagni}: The former bypasses complicated, often impossible derivation of stochastic representative solutions required for the latter, and, in addition, provides an approximate input-response (GF) function, applicable to generic IBVP's, subject to any combination of Dirichlet, Neumann, and/or Robin BC's. A simple naive estimator \cite{silverman, milsteindensity, devroye1985nonparametric, olariu, rudemo} for the 
transition density, $ p \left( \mathbf{x}, t| \mathrm{x'}, t' \right) , $ is then derived and heuristically validated using a central limit argument.
Numerical validation tests reveal two significant sources of error: depletion of the initial random walker swarm due to absorption at Dirichlet boundaries, and spatial discontinuity in estimated Green's function's, produced by our naive density estimator. To solve the first problem, we propose a simple particle splitting technique which preserves the number of simulated random walker's launched, as well as satisfying continuum mass conservation at the splitting point.  A simple area-averaging technique is proposed for tackling the second problem.    
In order to stochastically estimate GF's, three numerical parameters must be chosen: the number of random walker's launched, $ N_{\mathbf{x}, t} ,$ from a chosen response point, the interrogation area, $ \Delta A' \left( \mathbf{x'} , t' \right) , $ on which random walker's are captured, and the time step, $ \Delta t' , $ used for estimating random walker trajectories. Strategies for choosing these parameters are outlined in S1 Appendix A. Finally, the Green's function estimation technique is tested against three test problems, of increasing complexity, each having known GF's.    

\section{Recipe: Stochastic estimation of Green's functions}
Using stochastic methods to estimate GF's requires a number of analytical and numerical steps. In order to provide an overview of the proposed technique, we outline the steps here. 

\vspace{0.3cm}

\noindent \textbf{Step 1: Simultaneously determine the adjoint problem governing the Green's function and the magic rule:} Using Green's identities, and by trial and error, simultaneously determine the adjoint equation, the magic rule \cite{cole, barton1989elements}, and the boundary and initial conditions on the adjoint equation.  We illustrate this step in Appendix C, using a simplified, constant coefficient advection-diffusion problem.

\vspace{0.3cm}

\noindent \textbf{Step 2: Postulate that continuum transport reflects microscale stochastic dynamics:} This step provides a physical connection between microscale and continuum transport, as well as the mathematical basis for Steps 4 and 5. Express $ \eta \left( \mathbf{x'}, t' \right) $ as a functional, $ \eta \left( \boldsymbol{\chi} \left( s' \right) , s' \right) ,$ of a Brownian stochastic process,  $ \boldsymbol{\chi} \left( s' \right) $ \cite{gardiner2004handbook, schuss1, schuss2},
\begin{equation}\label{etachiassump}
\eta = \eta \left( \mathbf{\chi} \left( s' \right) , s' \right)
\end{equation}
evolving in backward time, $ s' , $ according to
\begin{equation}\label{stocheqn}
d \boldsymbol{\chi} \left( s' \right) = \mathbf{b} \left( \boldsymbol{\chi} \left( s' \right) ,s' \right) ds' + \mathbf{B} \left( \boldsymbol{\chi} \left( s' \right),s' \right) \cdot d \mathbf{w} \left( s' \right)
\end{equation}
where $ \mathbf{w} \left( s' \right) $ is an n-dimensional Weiner process, $ \mathbf{b} \left( \mathbf{x'} , s' \right) $ is the reversed drift/velocity field, $ \mathbf{b} \left( \mathbf{x'} , s' \right) = - \mathbf{v} \left( \mathbf{x'} , s' \right) ,$ and where the diffusion matrix is tied to $ \mathbf{B} $ by 
\begin{equation}\label{DandB}
\mathbf{D} \left( \mathbf{x'} , s' \right) = \langle \mathbf{B} \left(  \boldsymbol{\chi} \left( s' \right) , s' \right) \cdot \mathbf{B^T} \left( \boldsymbol{\chi} \left( s' \right) , s' \right) \rangle / 2
\end{equation}
with the average representing the conditional expectation $ \mathrm{E_{\mathbf{x}, s}} \mathbf{B} \left( \boldsymbol{\chi} \left( s' \right),s' \right) \cdot \mathbf{B^T} \left( \boldsymbol{\chi} \left( s' \right),s' \right) .$  In addition, the relationships between backward and forward (actual) time instants are as follows: i) the chosen, fixed (actual) response time, $ t , $ corresponds to the fixed backward instant, $ s ,  $ from which the random walker swarm is launched: $ t \longleftrightarrow s ;$ ii) the variable (actual) impulse time, $ t',$ corresponds to the variable backward random walker capture time, $ s' :$ $ t' \longleftrightarrow s' ;$ iii) the magnitude of forward time difference, $ t - t' \geq 0 ,$ equals the backward time difference, $ s' - s \geq 0 . $     

\vspace{0.3cm}

\noindent \textbf{Step 3: Determine the backward Kolmogorov problem governing backward evolution of the transition density:}   
The probability of observing a backward moving random walker, $ \boldsymbol{\chi} \left( s' \right) ,$ governed by (\ref{stocheqn}), at $ \left( \mathbf{x'},s' \right) ,$ given that it started at $ \left( \mathbf{x},s \right) , $ is given by the probability density function, $ p \left( \mathbf{x}, s | \mathbf{x'} , s' \right) .$ As discussed in Step 4, the present step is required in order to ensure that the adjoint problem governing the Green's function coincides with the terminal-boundary value problem governing $ p \left( \mathbf{x}, s | \mathbf{x'} , s' \right) .$ 

The backward time evolution of $ p \left( \mathbf{x}, s | \mathbf{x'} , s' \right) $ is governed by the backward Kolmogorov equation \cite{gardiner2004handbook, schuss1, schuss2}:
\begin{equation}\label{backwardkolmogorov}
\frac{\partial p \left( \mathbf{x} , s | \mathbf{x'},s' \right)}{\partial s'} + \mathbf{b} \left( \mathbf{x'},s' \right) \cdot \nabla' p \left( \mathbf{x} , s | \mathbf{x'},s' \right) + D_{ij} \frac{ \partial^2 p \left( \mathbf{x} , s | \mathbf{x'},s' \right)}{\partial x_i' \partial x_j'} =0
\end{equation}
where $ p \left( \mathbf{x} , s | \mathbf{x'},s' \right) $ is subject to the terminal backward time condition:
\begin{equation}\label{sametime}
p \left( \mathbf{x} , s | \mathbf{x'},s' \right) \rightarrow \delta \left( \mathbf{x} - \mathbf{x'} \right) \hspace{0.7cm} s' \rightarrow s
\end{equation}

Boundary conditions on (\ref{backwardkolmogorov}) that are relevant to, and correspond respectively to Dirichlet, Neumann, and Robin BC's in typical advection-diffusion-reaction problems, are as follows \cite{schuss1, schuss2}:
\begin{equation}\label{backwardkolmogorovdirichletbc}
p \left( \mathbf{x} , s | \mathbf{x'},s' \right) =0 \hspace{0.3cm} \mathrm{on} \hspace{0.3cm} \delta \Omega_D \times [s,T) 
\end{equation}
\begin{equation}\label{backwardkolmogorovneuwmannbc}
\frac{ \partial p \left( \mathbf{x} , s | \mathbf{x'},s' \right)}{\partial \tilde{n} \left( \mathbf{x'}, s \right)} =0 \hspace{0.3cm} \mathrm{on} \hspace{0.3cm} \delta \Omega_N \times [s,T) 
\end{equation}
and
\begin{equation}\label{backwardkolmogorovrobinbc}
- \frac{ \partial p \left( \mathbf{x} , s | \mathbf{x'},s' \right)}{\partial \tilde{n} \left( \mathbf{x'}, s \right)} = \kappa \left( \mathbf{x'},s' \right) p \left( \mathbf{x} , s | \mathbf{x'},s' \right)
\hspace{0.3cm} \mathrm{on} \hspace{0.3cm} \delta \Omega_R \times [s,T) 
\end{equation}
In Eq. (\ref{backwardkolmogorovdirichletbc}), $ \delta \Omega_D \times [s,T) $ denotes the space-time boundary on which Dirichlet BC's are imposed on the physical transport problem.  Similar definitions for Neumann and Robin boundaries are given respectively in (\ref{backwardkolmogorovneuwmannbc}) and (\ref{backwardkolmogorovrobinbc}). The conormal vector is defined as \cite{schuss2} $ \tilde{n} \left( \mathbf{x'}, s \right) = \mathbf{D} \left( \mathbf{x'}, s \right) \cdot \mathbf{n} \left( \mathbf{x'} \right),$ where $ \mathbf{n} \left( \mathbf{x'} \right) $ is the local unit outward normal.  Finally, in mass transfer problems, $ \kappa \left( \mathbf{x'},s' \right) $ is the local reaction coefficient \cite{schuss2}, arising, for example, on reactive boundaries. In heat transfer problems, $ \kappa \left( \mathbf{x'},s' \right) $ corresponds to the local convective heat transfer coefficient, divided by the local thermal conductivity.

\vspace{0.3cm}

\noindent \textbf{Step 4: Determine the applicability of the Green's function estimation method:} For transport problems in which diffusive or dispersive transport is non-negligible, stochastic estimation of GF's can only be used if the adjoint problem governing the so-called propagator, $ K \left( \mathbf{x} , s | \mathbf{x'},s' \right), $ coincides with the backward Kolmogorov problem, governed by Eq. (\ref{backwardkolmogorov}) and problem-specific combinations of the BC's in Eqs. (\ref{backwardkolmogorovdirichletbc}) through (\ref{backwardkolmogorovrobinbc}). Specifically: 

\vspace{0.3cm}

\noindent a) the diffusion/dispersion matrix must be diagonal:
\begin{equation}\label{diffntensor}
\mathbf{D} = \mathbf{\tilde{D}} \left( \mathbf{x} , s \right) \mathbf{\cdot} \mathbf{I}    
\end{equation}
This requirement emerges in Step 1: Non-diagonal $ \mathbf{D} $ produces a convolution integral, over the problem domain, $ \Omega , $ in two unknowns, the Green's function and $ \eta \left( \mathbf{x'} , s' \right) . $

\vspace{0.3cm}

\noindent b) The spatial domain, $ \Omega , $ in which transport takes place, must remain fixed/non-deforming; this ensures that time-invariance and symmetry properties intrinsic to the propagator
are enforced \cite{barton1989elements}.

Under conditions a) and b), the problem governing the propagator, $ K \left( \mathbf{x} , s | \mathbf{x'},s' \right) ,$  mirrors the problem governing $ p \left( \mathbf{x} , s | \mathbf{x'},s' \right) ,$ with $ K \left( \mathbf{x} , s | \mathbf{x'},s' \right) $ replacing $ p \left( \mathbf{x} , s | \mathbf{x'},s' \right) $ in Eqs. (\ref{backwardkolmogorov}) through (\ref{backwardkolmogorovrobinbc}).  Given an estimated propagator, the associated Green's function is then given by \cite{barton1989elements}:
\begin{equation}\label{propandgf}
G \left( \mathbf{x} , s | \mathbf{x'},s' \right) = H \left( s' - s \right) K \left( \mathbf{x} , s | \mathbf{x'},s' \right)
\end{equation}
where $H \left( s' - s \right) $ is the Heaviside function.

\vspace{0.3cm}

\noindent \textbf{Step 5: Choose a scheme for integrating Eq. (\ref{stocheqn}) and an estimator for estimating $ p \left( \mathbf{x}, s | \mathbf{x'} , s' \right) :$}  Various approaches, each characterized by a pointwise (strong) order of convergence, are available for numerically integrating stochastic differential equations, here  Eq. (\ref{stocheqn}) \cite{milstein1, milstein2, kloeden, sdeint1}.  Since boundary conditions on the adjoint equation are all-determining in constructing the Green's function, and since interaction of stochastic processes with Dirichlet, Neumann, and Robin boundaries, as integrated by the Euler-Maruyama scheme \cite{milstein1, kloeden, sauer} have  been well-studied \cite{schuss2}, this scheme is recommended.  We note below, however, alternative, less expensive integration schemes.

Estimation of probability density functions, here $ p \left( \mathbf{x}, s | \mathbf{x'} , s' \right) ,$ is a well-developed field \cite{silverman, milsteindensity, devroye1985nonparametric, olariu, rudemo}. In this paper, as detailed below, we use a simple naive estimator \cite{silverman, devroye1985nonparametric}.  Limitations and alternatives to this approach are highlighted below.

\section{Comparison of proposed Green's function estimation method with Monte Carlo solutions of initial-boundary value problem; limitations}
Monte Carlo (MC) techniques for solving linear partial differential equations, subject to problem-specific boundary and initial conditions, have a long history \cite{ulam, metropolis, todd}.  By contrast, stochastic construction of GF's via transition density estimation appears to represent a new and more general approach, where the estimated Green's function can be used to (approximately) solve a given linear PDE, subject to any combination of initial and boundary conditions.  In order to contrast the proposed stochastic Green's function estimation method against MC solutions of IBVP's \cite{muller,sabelfeld,booth,sabelfeldsimonov,haji,dimnov}, we outline the five-step construction of MC solutions, and then touch on the significant differences and advantages introduced by the present approach. 

Monte Carlo solutions are based on representative stochastic solutions of elliptic, parabolic, and hyperbolic PDE's, subject to specific boundary conditions, and in the latter two cases, initial conditions \cite{feynmanphd, kac1, kac, friedman1975stochastic, bass}. Focusing on parabolic problems  governed, for example, by Eqs. (\ref{mainequation1}) or (\ref{mainequation}), the recipe is as follows: 
i) Assume that macroscopic evolution of the field variable, $ \eta \left( \mathbf{x} , t \right) ,$  reflects (long time scale, collective \cite{ keanini2011green, schuss2}) microscale advective-diffusive transport of some conserved real or virtual property, e.g., microscopic number, momentum, or energy density  \cite{keanini2011green, forster, boonyip}, represented by $ \boldsymbol{\chi} \left( s \right) ,$ and modeled by Eq. (\ref{stocheqn}). ii) Using Ito's formula \cite{gardiner2004handbook, schuss1, schuss2}, determine the differential, backward-time change in $ \eta , $ produced by the stochastic evolution of $ \boldsymbol{\chi} \left( s \right). $ iii) In the equation derived in step ii), and referring to (\ref{mainequation1}), replace the temporal derivative of $ \eta , $ along with advective and diffusive terms with remaining source and reaction terms, $ f \left( \mathbf{x}, t \right) $ and $ \gamma \eta \left( \mathbf{x}, t \right) ,$ respectively. iv) Attempt to obtain a stochastic representative solution for $ \eta \left( \mathbf{x}, t \right) , $ by integrating, backward in time, the differential equation obtained in iii). In deriving the representative solution, account for the stochastic process sampling Dirichlet, Neumann, and/or Robin boundaries that appear in the continuum transport problem \cite{kac1, kac, friedman1975stochastic, bass}. v) Using the representative solution obtained in step iv), launch a random walker swarm from a chosen solution point, $ \left( \mathbf{x}, t \right) ,$ and use appropriate absorption (Dirichlet boundaries), reflection (Neumann boundaries), and partial reflection (Robin boundaries) techniques \cite{schuss2} to construct a Monte Carlo estimate for $ \eta \left( \mathbf{x}, t \right). $ 

Comparing the MC approach with the proposed Green's function estimation procedure, two significant advantages emerge. First, by directly estimating transition densities $ p \left( \mathbf{x}, s | \mathbf{x'} , s' \right) ,$ Green's function estimation bypasses the challenging task of deriving representative stochastic solutions \cite{kac1, kac, friedman1975stochastic, bass}: Steps iii) and iv) are eliminated in favor of straightforward stochastic integration of Eq. (\ref{stocheqn}), combined with well-developed boundary interaction techniques \cite{schuss2}.
Second, as noted above, the technique produces an estimated, non-problem-specific input-response function (Green's function) that can be used to solve a family of advection-diffusion-reaction problems, subject to any combination of Dirichlet, Neumann, and/or Robin BC's.

\subsection{Challenges and a limitation}
The present work is experimental, focused on developing a stochastically-based Green's function estimation technique. Prior to looking at details, we highlight two significant challenges that required solution, and note a limitation associated with the method in its current form.  

First, following launch of a random walker swarm - consisting of $ N_{\mathbf{x}, s} $ integrated realizations of Eq. (\ref{stocheqn}) - from a chosen solution point, $ \left( \mathbf{x} , s \right) , $ absorption of individual realizations at Dirichlet boundaries rapidly depletes the initial set of $ N_{\mathbf{x}, s} $ random  walkers. Depletion, in turn, produces increasingly degraded-in-time Green's function estimates. Our solution, based on respawning of weighted random walker's, and carried out at random locations in the problem domain, enforces continuum conservation of mass at the particle respawning point, and significantly improves Green's function estimation accuracy. See Section \ref{respawnsection}. Second, spatial variance in estimated GF's likely arises due to use of a spatially discontinuous naive density estimator \cite{silverman, devroye1985nonparametric}, Eq. (\ref{naiveestp}) below. We address this problem by introducing an area-averaging technique in which variance in local Green's function estimates is minimized by iteratively altering the size of a smoothing window. See Section \ref{areavgsection}. Although not investigated, kernel estimators, which also eliminate discontinuity in naive estimators \cite{silverman, devroye1985nonparametric}, might also prove advantageous.

Regarding the present technique's main limitation, Euler-Murayama integration of Eq. (\ref{stocheqn}) has a low $ \Delta s ^{1/2} $ order of strong convergence \cite{milstein1, milstein2, kloeden}.  For a given degree of Green's function estimation accuracy, significantly faster integration can be achieved using, for example, Milstein's order 1 method \cite{milstein2}, or order 2 strong Taylor expansion methods \cite{kloeden}, or order 2 Runge-Kutta methods \cite{sdeint1}. However, again, since boundary conditions on the adjoint problem are central to constructing the estimated Green's function, and since treatment of random walker-boundary interactions has reached maturity in the case of Euler-Murayama integration \cite{schuss2}, this scheme is used.

\section{A simple estimator for the transition density}\label{transitiondensityest}
Transition density estimation remains a research-level task \cite{silverman,milsteindensity,devroye1985nonparametric,olariu, rudemo}. An accessible overview of theoretical results, focused on the simplest case of estimating 'free-space' transition densities, $ p_{o} \left( \mathbf{x} , s | \mathbf{x'},s' \right) , $ can be found in \cite{milsteindensity}. 



In this study, the transition density is estimated using the naive estimator \cite{silverman}:
\begin{equation}\label{naiveestp}
p \left( \mathbf{x} , s | \mathbf{x'},s' \right) = \frac{n \left( \mathbf{x'}, s' \right)}{N_{\mathbf{x},s} \Delta A'\left( \mathbf{x'}, s' \right)} +  O \left( {\Delta s'} \right)
\end{equation}
where, the $  O \left( {\Delta s'} \right) $ estimation error is determined by the Euler-Murayama scheme's weak convergence order \cite{gardiner2004handbook, milstein2, kloeden}, and where $ N_{\mathbf{x},s} ,$ $ \Delta A'\left( \mathbf{x'}, s' \right), $ and $ n \left( \mathbf{x'}, s' \right) $ are, respectively, the number of RW's launched from the solution point, $ \left( \mathbf{x}, s \right) ,$ the chosen RW interrogation area, centered on $ \mathbf{x'} ,$ at backward time, $ s' ,$ and the number of RW's that reach or pass through $ \Delta A'\left( \mathbf{x'}, s' \right), $ over the time interval, $ [ s' , s' + \Delta s' ) . $ 


While the formula in Eq. (\ref{naiveestp}) is stated in terms of the Lebesque measure in, e.g., \cite{devroye1985nonparametric}, it can be derived heuristically as follows.  Focus on two-dimensional random walks, as used in our tests, evolving within an arbitrary two-dimensional domain, $ \Omega , $ subject to any combination of Dirichlet, Neumann and/or Robin boundaries. For two-dimensional problems, it is useful to picture a three-dimensional, cylindrical hyperspace in which the downward (upward) vertical direction corresponds to the backward (forward) time axis, $ s' \left( \leftrightarrow t' \right) ,$ the top corresponds to the backward (forward) solution time slice, $ s \left( \leftrightarrow t \right), $ the bottom to the backward (forward) terminal (initial) time slice, with any horizontal slice through the hyperspace corresponding to an instantaneous snapshot of the two-dimensional spatial domain, $ \Omega .$ 


Symbolically, recognize that the probability of random walker's, $ \boldsymbol{\chi} \left( s' \right) , $ launched from $ \left( \mathbf{x}, s \right) , $ reaching a chosen interrogation area, $ \Delta A' \left( \mathbf{x'}, s' \right) , $ at the instant $s' , $ is given by:
\begin{equation}\label{probinterrogation}
P \left[ \boldsymbol{\chi} \left( s' \right) \in \Delta A' \left( \mathbf{x'}, s' \right) \right] = \int_{\Delta A' \left( \mathbf{x'}, s' \right) } p \left( \mathbf{x} , s | \mathbf{x"},s' \right) d \mathbf{x"} 
\end{equation}
Assuming that $ p \left( \mathbf{x} , s | \mathbf{x'},s' \right) $ is well-behaved, then in the limit as $ \Delta A' \left( \mathbf{x'}, s' \right) \rightarrow 0 , $ by the mean-value theorem,
\begin{equation}\label{limitp}
\int_{\Delta A' \left( \mathbf{x'}, s' \right) } p \left( \mathbf{x} , s | \mathbf{x"},s' \right) d \mathbf{x"} = p \left( \mathbf{x} , s | \mathbf{x^*},s' \right) \Delta A' \left( \mathbf{x'}, s' \right) 
\end{equation}
where $ \mathbf{x^*} \in \Delta A' \left( \mathbf{x'}, s' \right) ,$ and $ \mathbf{x^*} \rightarrow \mathbf{x'} $ as $ \Delta A' \left( \mathbf{x'}, s' \right) \rightarrow 0 . $

Thus, based on the numerical simulation of random walker evolution following launch, approximate the probability on the left side of Eq. (\ref{probinterrogation}) as
\begin{equation}\label{probareaapprox}
P \left[ \boldsymbol{\chi} \left( s' \right) \in \Delta A' \left( \mathbf{x'}, s' \right) \right] \approx \frac{n \left( \mathbf{x'}, s' \right)}{N_{\mathbf{x},s}}
\end{equation}
Combining Eqs. (\ref{probinterrogation}) through (\ref{probareaapprox}) then leads to the estimator in (\ref{naiveestp}). As highlighted in \cite{silverman}, and as we find below, estimators like Eq. (\ref{naiveestp}) produce discontinuous estimates of $ p \left( \mathbf{x} , s | \mathbf{x'},s' \right) . $ 

\subsection{Heuristic validation of estimator, Eq. (\ref{naiveestp})}\label{centrallimit}
We adapt a central limit argument in \cite{kubobook2} to estimate the free-space transition density, $ p_{o} \left( \mathbf{x} , s | \mathbf{x'},s' \right) , $ for the evolution of a numerically simulated swarm of $ N_{\mathbf{x},s} $ RW's, launched from $ \left( \mathbf{x}, s \right) , $ within an infinite, two-dimensional domain.  Of the $ N_{\mathbf{x},s} $ RW's launched, focus on the $ M $ RW's that reach or pass through differential area, $ \Delta A' \left( \mathbf{x'}, s'\right) ,$ over the backward time interval, $ n \Delta s' \leq s' < \left( n+1 \right) \Delta s',$ where $ \mathbf{x'} \in \Delta A' \left( \mathbf{x'}, s'\right) ,$ and $ s' = s + n \Delta s .$

Define the space-time vector displacement of the $ j^{th} $ random walker, in the set of $ M $ RW's captured by $ \Delta A' \left( \mathbf{x'}, s'\right) , $ as 
\begin{equation}\label{vectorg}
\mathbf{g_n^{\left(j\right)}} \left( s' \right) =  \sum_{k=1}^n \boldsymbol{\Delta} \boldsymbol{\chi_k^{\left(j\right)}} = \boldsymbol{\chi^{\left(j\right)}} \left( s' \right) - \mathbf{x}  \hspace{1cm} j = 1, 2, ... M
\end{equation}
or, in component form,
\begin{align}\label{gdefnx}
g_{n_x}^{\left( j \right)} & = & \Delta x_1^{\left( j \right)} + \Delta x_2^{\left( j \right)} + ... + \Delta x_n^{\left( j \right)} \hspace{4cm} \nonumber \\
  & = & \chi_x^{\left(j\right)} \left( s' \right) - x \hspace{6.7cm}  \nonumber \\
g_{n_y}^{\left( j \right)} & = & \Delta y_1^{\left( j \right)} + \Delta y_2^{\left( j \right)} + ... + \Delta y_n^{\left( j \right)}  \hspace{4cm} \nonumber \\
 & = &  \chi_y^{\left(j\right)} \left( s' \right) - y  \hspace{6.7cm} 
\end{align}
where 
\begin{align}
\boldsymbol{\Delta} \boldsymbol{\chi_k^{(j)}}  & = & \boldsymbol{\chi^{(j)}} \left( \left(k+1\right) \Delta s' \right) -  \boldsymbol{\chi^{(j)}} \left( k \Delta s' \right)  \hspace{9cm} \nonumber \\
 & = & \sqrt{2 D_o } \left[ \mathbf{W^{(j)}} \left( \left(k+1\right) \Delta s' \right) - \mathbf{W^{(j)}} \left( k \Delta s' \right)  \right] \hspace{7cm} 
\end{align}
is the simulated Brownian displacement of random walker $ j ,$ over the $ k^{th} $ time-interval.

Now, let $ N_{\mathbf{x}, s} \rightarrow \infty , $ so that, likewise, $ M \rightarrow \infty , $ and define the random variable $ \tilde{Y}_{n_x} $ as the sum of $ n $ independent random variables, $ \Delta x_1 ,$ $ \Delta x_2 , ..., \Delta x_n :$
\begin{equation}\label{ydefn}
\tilde{Y}_{n_x} = \frac{g_{n_x}}{s_{n_x}}= \frac{\Delta x_1 + \Delta x_2 + ... + \Delta x_n}{s_{n_x}}
\end{equation}
where $ s_{n_x}^2 = \langle \Delta x_1^2 \rangle +\langle \Delta x_2^2 \rangle + ... + \langle \Delta x_1^n \rangle ,$ is the sum of variances of the set of $ n $ differential x-displacements.  Importantly, recognize that the differential spread, i.e., variance, in backward time x-displacements (over the set of $ M $ random walker trajectories captured on $ \Delta A' \left( \mathbf{x'}, s' \right) )$ \textit{decreases} with increasing backward time: $ \langle \Delta x_1^2 \rangle = -2D_o \Delta s' , $ or, in terms of forward time increments, $ \langle \Delta x_1^2 \rangle = 2D_o \Delta t' , $ 

Letting $ \Delta s'  \rightarrow 0 ,$ so that $ n \rightarrow \infty , $ then by the central limit theorem \cite{kubobook2}, the probability density for $ g_{n_x} $ tends to the Gaussian:
\begin{equation}\label{pdfgx}
P \left( g_{n_x} \right) \rightarrow \frac{1}{\sqrt{4 \pi D_o \left( t - t' \right)}} \exp{ - \frac{ \left( x - \langle \chi_x \left( t' \right) \rangle \right)^2}{4 D_o \left( t - t'\right) }}
\end{equation}
where, since $ \langle \Delta x_1^k \rangle = 2 D_o \Delta t' ,$ $ k =1,2,..., n ,$ $ s_{n_x}^2 = 2D_o \left( t -t'  \right) .$ The same argument applied to the y-component of the random relative displacement, $ g_{n_y} $ leads to:
\begin{equation}\label{pdfgy}
P \left( g_{n_y} \right) \rightarrow \frac{1}{\sqrt{4 \pi D_o \left( t - t' \right)}} \exp{ - \frac{ \left( y - \langle \chi_y \left( t' \right) \rangle \right)^2}{4 D_o \left( t - t'\right) }}
\end{equation}
where again,  $ t \leftrightarrow s , $ $ t' \leftrightarrow s' , $ and $ t \geq t ' . $  

Finally, letting $ \Delta A' \left( \mathbf{x'}, s' \right) \rightarrow 0,$ then $ \langle \chi_x \left( s' \right) \rangle \rightarrow x' $ and $ \langle \chi_y \left( s' \right) \rangle \rightarrow y' .$  Thus, we find that the simple estimator in Eq. (\ref{naiveestp}) for the probability of observing a random walker at $ \left( \mathbf{x'}, s' \right) ,$ given that it was launched at $ \left( \mathbf{x}, s \right) ,$ given by the product of the right sides of Eqs. (\ref{pdfgx}) and (\ref{pdfgy}), corresponds to the solution of the diffusive backward Kolmogorov equation, Eq. (\ref{backwardkolmogorov}): 
\begin{equation}
p \left( \mathbf{x} , t | \mathbf{x'} , t' \right) = \frac{1}{\sqrt{4 \pi D_o \left( t - t' \right)}} \exp{ - \frac{ \left( \mathbf{x} - \mathbf{x'} \right)^2}{4 D_o \left( t - t'\right) }}
\end{equation}
Equivalently, and as required \cite{gardiner2004handbook, schuss2}, this corresponds to the free-space solution to the (pure diffusion) Fokker-Planck equation, governing forward time evolution of $ p \left( \mathbf{x} , t | \mathbf{x'} , t' \right),$ following delta function release of RW's, at $ t= t',$ from $ \mathbf{x} = \mathbf{x'} .$ 

\subsection{Green's function estimate and error in the estimate}
Assuming that the adjoint problem governing the propagator, $ K \left( \mathbf{x} , s | \mathbf{x'},s' \right),$ and the backward Kolmogorov problem governing the transition density have the same structure, then from Eq. (\ref{propandgf}), for $ s' - s \geq 0 , $ the stochastically estimated Green's function follows immediately from Eq. (\ref{naiveestp}):
\begin{equation}\label{stochestgf}
G \left( \mathbf{x} , s | \mathbf{x'},s' \right) = p \left( \mathbf{x} , s | \mathbf{x'},s' \right) = \frac{n \left( \mathbf{x'}, s' \right)}{N_{\mathbf{x},s} \Delta A'\left( \mathbf{x'}, s' \right)} +  O \left( {\Delta s'} \right)
\end{equation}
where the Euler-Murayama scheme's $  O \left( {\Delta s'} \right) $ weak convergence is again highlighted.

We note that an additional, localized, $ O \left( \Delta s' \right) $ error arises in connection with the random walker respawning algorithm introduced below.  As shown in Appendix B,  this is an implicit discretization error, confined to random locations in the solution domain where individual RW's are respawned following absorption on Dirichlet boundaries.  

\section{Absorbing boundaries and random walker respawning}
In order to mitigate against loss of absorbed RW's and preserve early-forward-time solution accuracy, we introduce a random walker respawning algorithm \cite{luxsplit}, and discuss below and in Appendix B its properties:

\vspace{0.4cm}

\noindent a) Initially, at the chosen backward solution instant, $ s, $ the group of $ N_{\mathbf{x}, s} , $ RW's to be launched from the solution point, $ \mathbf{x} , $ are all assigned a weight of 1, and each random walker weight is stored in an $ N_{\mathbf{x},s} \times 1 -\mathrm{dimensional} $ array, $ [ W ] $ \cite{luxsplit}.

\vspace{0.4cm}

\noindent b) After every backward time-step, $ \Delta s' , $ any random walker that reaches a Dirichlet boundary is removed, and the in-domain random walker, $ \chi \left( s' ; W_{max} \right) , $ having the highest weight, $ W_{max} \left( s' \right) , $ in $ [ W ] $ is also removed and replaced by two RW's, each having weights equal to $ W_{max} \left( s' \right) / 2 . $ The replacement pair are placed at the same location occupied by $ \chi \left( s' ; W_{max} \right) .$  

\vspace{0.4cm}

\noindent c) After each set of $ m $ time-steps, the array $ [ W ] $ is reordered, in descending order, according to the current set of random walker weights. Thus, high-mass RW's that experience relatively few absorption events, in particular, those far from Dirichlet boundaries, tend towards Dirichlet boundaries over backward time, suppressing numerical thinning of near-boundary mass distributions. 

\vspace{0.4cm}

When using the respawning algorithm, equation (\ref{stochestgf}) must be modified:
\begin{equation}\label{gengfstochasticconstructionrespawn}
G\left( \mathbf{x}, t | \mathbf{x}^{\prime} , t^{\prime} \right) = G\left( \mathbf{x}, s | \mathbf{x}^{\prime} , s^{\prime} \right) = \frac{1}{N_{\mathbf{x} , s}} \frac{1}{ \Delta A^{\prime} } W \left( \mathbf{x}^{\prime} , s^{\prime} \right) + O \left( {\Delta s'} \right)
\end{equation}
where $W(\mathbf{x}',s')$ is the total weight (sum of weights) of the RW's that reach $\Delta A' (\mathbf{x}',s')$. 

\subsection{Properties of the respawning algorithm}
Appendix B shows that the respawning algorithm has two important properties which ensure the physical consistency of the algorithm, as well as a third that significantly enhances the statistical quality of the Green's function estimation procedure:

\vspace{0.4cm}

\noindent a) The algorithm enforces, at each random respawning location, the continuum diffusion equation governing continuum random walker evolution.

\vspace{0.4cm}

\noindent b) Since the algorithm diffusively smooths, on short, single-time-step time-scales, instantaneous, localized mass gradients produced by respawning, the algorithm appears to be insensitive to the number of time steps, $ m ,$ taken between reordering of the random walker weight array, $ [ W ] .$  In this first study, we have not performed tests to ascertain suitable ranges for $ m ;$ in all cases, $ m = 10 .$

\vspace{0.4cm}

\noindent c) Crucially, the respawning algorithm conserves the number of random walkers, $ N_{\mathbf{x}, s} , $ initially launched.  As shown in the Results below, preserving $ N_{\mathbf{x}, s} $ is essential to obtaining low-error Green's function estimates.

\section{Numerical experiments}
  
We focus on three test cases, each of increasing complexity: 1) two-dimensional diffusive transport in a square, 2) diffusion in a circle, and 3) unsteady two-dimensional groundwater solute transport, produced by combined advection,  dispersion, and chemical depletion in a nonhomogeneous, anisotropic, infinite region. We use the first test to expose and solve two principle challenges: a) loss of random walker's at Dirichlet boundaries and b) spatial variance in estimated GF's. The second test examines the performance of the Green's function estimation procedure in domains having non-planar boundaries.  Since tests 1 and 2 allow (preferred, optimal) backward time Green's function construction, in the next section, we illustrate, in detail, the backward time Green's function estimation recipe.  The third test requires forward time Green's function estimation; there, we briefly outline the corresponding recipe. 


\subsection{Application of the backward time Green's function estimation recipe to pure diffusion test cases 1 and 2}

\subsubsection{Step 1: Simultaneously derive the adjoint equation and magic rule}
Appendix C illustrates trial and error derivation of these two central relations, focusing on a generic advection-diffusion transport problem in a finite domain.  Similar trial and error approaches were used to derive the adjoint equations and magic rules used in all three test cases. In tests 1 and 2, the adjoint equation is thus given by:
\begin{equation}\label{adjointeqnnew}
\frac{\partial G \left( \mathbf{x} , s | \mathbf{x'},s' \right)}{\partial s'} + D_{o} {\nabla'}^2 G \left( \mathbf{x} , s | \mathbf{x'},s' \right) = - \delta \left( \mathbf{x} - \mathbf{x'} \right)  \delta \left( s -s' \right) 
\end{equation}
where $ D_o $ is a fixed diffusion coefficient. 

Considering the case where a combination of Dirichlet and Neumann boundary conditions are imposed, the magic rule assumes the form:
\begin{multline}\label{magicrulenew}
\eta(\mathbf{x},t) = \int_0^t \int_{\Omega} G(\mathbf{x},t|\mathbf{x'},t')f(\mathbf{x'}, t)d \mathbf{x'}dt' - D_o \int_0^t \int_{\partial \Omega_D} g(\mathbf{x'},t') \nabla' G(\mathbf{x},t|\mathbf{x'},t') \cdot \mathbf{\hat{n}'} d\mathbf{x'}dt' + \\
+ D_o \int_0^t \int_{\partial \Omega_N} G(\mathbf{x},t|\mathbf{x'},t') \nabla' h(\mathbf{x'},t') \cdot \mathbf{\hat{n}'}  d\mathbf{x'}dt' + \\
+ \int_{\Omega} \phi(\mathbf{x'})G(\mathbf{x},t|\mathbf{x'},0)d \mathbf{x'}.
\end{multline}
where $ f(\mathbf{x'}, t),$ $ g(\mathbf{x'},t'), $ $ h(\mathbf{x'},t'),$ and $ \phi(\mathbf{x'})$ represent, respectively, the areal/volumetric source/sink term, boundary conditions on the Neumann $ ( \partial \Omega_N ) $ and Dirichlet $ ( \partial \Omega_D ) $ boundaries, and the initial condition. See Eq. (\ref{magicsuppinfo2}), Appendix C, with the advective (last) term turned off. In addition, $ \mathbf{\hat{n}'}$ is the local outward unit normal to the problem domain, $ \Omega' .$  
\subsubsection{Steps 2 and 3: Postulate existence of a microscale stochastic transport process, $ \boldsymbol{\chi} \left( s \right), $ and write down the  associated backward Kolmogorov problem}
Since test cases 1 and 2 model pure diffusion, the
backward-time stochastic differential equation is obtained by dropping the drift term in Eq. (\ref{stocheqn}):
\begin{equation}\label{stocheqntest}
d \boldsymbol{\chi} \left( s \right) = \mathbf{B} \left( \boldsymbol{\chi} \left( s \right),s \right) \cdot d \mathbf{w} \left( s \right)
\end{equation}
where, from Eq. (\ref{DandB}), $ \langle \mathbf{B} \rangle = \sqrt{2 D_o} \mathbf{I} .$  Thus, the corresponding backward Kolmogorov equation, determined by dropping the drift term from Eq. (\ref{backwardkolmogorov}), is given by
\begin{equation}\label{backwardkolmogorovtest}
\frac{\partial p \left( \mathbf{x} , s | \mathbf{x'},s' \right)}{\partial s'} +  D_{o} {\nabla'}^2  p \left( \mathbf{x} , s | \mathbf{x'},s' \right) =0
\end{equation}
where $ p \left( \mathbf{x} , s | \mathbf{x'},s' \right)$ satisfies the terminal condition in Eq. (\ref{sametime}).
Finally, in tests 1 and 2, Dirichlet boundary conditions, as stated in Eq. (\ref{backwardkolmogorovdirichletbc}), are imposed at all points on both domain boundaries.  
\subsubsection{Step 4: Determine the applicability of the stochastic Green's function estimation procedure}
Comparing Eqs. (\ref{adjointeqnnew}) and (\ref{backwardkolmogorovtest}), due to the non-zero right side in (\ref{adjointeqnnew}), it is clear that these do not have the same mathematical
structure. However, the propagator, $ K \left( \mathbf{x} , s | \mathbf{x'},s' \right), $  which is connected to the Green's function, $ G \left( \mathbf{x} , s | \mathbf{x'},s' \right) $ via Eq.(\ref{propandgf}), 
\textit{does} satisfy the equation $ L^*K \left( \mathbf{x} , s | \mathbf{x'},s' \right)=0,$ where $ L^* = \partial / \partial s' + D_o {\nabla'}^2 = \partial / \partial t' + D_o {\nabla'}^2.$ Thus, using $ L^* $ to operate on $  H \left( \tau \right)  K \left( \mathbf{x} , s | \mathbf{x'},s' \right) ,$  where $ \tau= t -t' , $ using $ \partial H / \partial t' = - \delta \left( \tau \right) = - \delta \left( -\tau \right) = - \delta \left( s'- s \right) $ as well as the terminal condition, $  p \left( \mathbf{x} , s | \mathbf{x'},s' \right) \rightarrow \delta \left( \mathbf{x} - \mathbf{x'} \right) $ as $ s' \rightarrow t ,$ we obtain  Eq. (\ref{adjointeqnnew}).  

Thus, since the problems governing the propagator, $ K \left( \mathbf{x} , s | \mathbf{x'},s' \right) $ and transition density, $ p \left( \mathbf{x} , s | \mathbf{x'},s' \right) $ have identical mathematical structure, including matching BC's and IC's, we can use a density estimation procedure to first estimate $ K \left( \mathbf{x} , s | \mathbf{x'},s' \right) ,$ and then immediately obtain $ G \left( \mathbf{x} , s | \mathbf{x'},s' \right) $ via Eq. (\ref{propandgf}). 
\subsubsection{Step 5: Choose an integration scheme for the problem SDE and a density estimation method}
As noted, in all three tests, we use the Euler-Murayama scheme \cite{schuss2, milstein2, kloeden} to integrate Eq. (\ref{stocheqntest}).  Thus, discretize the backward time axis as: $ {s'}_{i+1} = {s'}_i + \Delta s' ,$ $ i=0, 1, 2, ... \mathrm{M} , $ where $ {s'}_o = s ,$ the chosen backward solution time, corresponds to a chosen forward solution time, $ t ,$ $s' $ is the chosen random walker interrogation time, and $ \mathrm{M} \Delta s' = s' - s . $ Defining $ \Delta \mathbf{w} \left( {s'}_i \right) = \mathbf{w} \left( {s'}_{i+1} \right) - \mathbf{w} \left( {s'}_i \right) ,$  as the psuedo-random, numerical approximation to the Weiner differential, $ d \mathbf{w} \left( {s'}_i \right) $ in Eq. (\ref{stocheqntest}), the EM scheme is given by \cite{schuss2, kloeden}:
\begin{equation}\label{eulerscheme}
\mathbf{y} \left( {s'}_{i+1} \right) = \mathbf{y} \left( {s'}_i \right) + \sqrt{ 2 D_o} \Delta \mathbf{w} \left( {s'}_i \right)
 \end{equation}
where the numerical approximation to the actual stochastic jump, $ d \boldsymbol{\chi} \left( {s'}_i \right) = \boldsymbol{\chi} \left( {s'}_{i+1} \right) - \boldsymbol{\chi} \left( {s'}_i \right) , $ in (\ref{stocheqntest}) is denoted as $\mathbf{y} \left( {s'}_{i+1} \right) - \mathbf{y} \left( {s'}_i \right) . $ Details on simulation of $ \Delta \mathbf{w} \left( {s'}_i \right) $ can be found, for example, in \cite{schuss2, kloeden}.



\subsection{Explanatory notes}

\noindent A) The stochastic estimation technique requires specification of three numerical parameters: i) the backward (forward) time-step size, $ \Delta s' $ $ \left( \Delta t \right) , $  ii) the interrogation area, $ \Delta A' \left( \mathbf{x}',s' \right) $ $ \left( \Delta A' \left( \mathbf{x}',s' \right) \right)  ,$ and iii) the number of RW's, $ N_{\mathbf{x},s} $ $ \left( N_{\mathbf{x'},t'}\right) , $ launched from the chosen Green's function solution point, $ \left( \mathbf{x}, s \right) $ $ \left( \left( \mathbf{x'}, t' \right) \right) .$ Appendix A proposes procedures for specifying these parameters. 


\vspace{0.4cm}

\noindent B) In all three tests, lengths, $ \left( x,y; x',y' \right), $ forward and backward times, $ \left( t, t' ; s, s' \right) , $ and diffusion coefficients, $ D_o ,$ are expressed in non-dimensional form: i) $ \left( x, y; x',y' \right) = \left( \tilde{x}, \tilde{y} ; \tilde{x}', \tilde{y}'  \right) \tilde{L}^{-1} , $ ii)$ \left( t,t'; s, s' \right) = \left( \tilde{t}, \tilde{t}'; \tilde{s}, \tilde{s}'  \right) \tilde{t}_s^{-1} , $ , iii) $ D_o = \tilde{D}_o \tilde{t}_s \tilde{L}^{-2} , $ where $ \tilde{L}  $ and $ \tilde{t}_s $ are dimensional, problem-specific length and time-scales. In the first two tests, $ \tilde{D}_o $ is a dimensional diffusion coefficient, and in the third, a dimensional dispersion coefficient.  Throughout the rest of this paper, all functions, variables, and parameters without a tilde are dimensionless.

In all three tests, we choose a dimensionless solution time, $ t =10 ,$ and dimensionless length, $ L = 1 . $ In order to gain a sense of associated dimensional scales, consider diffusion of small molecular and ionic species in water, where $ \tilde{D}_o = \mathrm{O} \left( 10^{-9} \ \mathrm{m^2} \mathrm{s^{-1}} \right); $ since diffusive spread of a specie introduced at a point grows as $ \tilde{L} \sim \sqrt{ \tilde{D}_o \tilde{t} } , $ then choosing time scales, $ \tilde{t}_s , $ equal to, e.g., $ 10^{-3} \ \mathrm{s} $ and  $ 10^{7} \ \mathrm{s} ,$ requires computational domains having dimensions on the order of $ \tilde{L} \sim \sqrt{ 10 \tilde{D}_o \tilde{t}_s } \sim 3 \left( 10^{-6} \right) \ \mathrm{m} $ and $ \sim 3 \left( 10^{-1} \right) \ \mathrm{m} ,$ respectively. Similar considerations can be used in the  groundwater transport problem in test 3.  




\vspace{0.4cm}

\noindent C) In test 1, we discretize the square solution domain, having dimensionless length, $ L=1 ,$ into $ 10^{4} $ (square) interrogation areas. In test 2, a circular solution domain of dimensionless diameter, $ D=1,$ is inscribed and centered within the same unit square. Again, the square is discretized into $ 10^{4} $ square interrogation areas.  Thus, in both tests, $ \Delta x' = \Delta y' = 10^{-2} ,$ and $ \Delta A_{\mathbf{x_j^{'}}, s'} = \Delta x' \Delta y' = 10^{-4} ,$ where $ \mathbf{x_j^{'}} $ denotes the centroid of the $ j^{th}$ interrogation area. While test 3 considers ground water transport in an infinite domain, for finite solution time, $ t = 10 , $ and under the groundwater flow conditions considered, we can again define a unit square solution domain, again having $ 10^4 $ interrogation areas; see Section \ref{test3section}.

In tests 1 and 2, the dimensionless diffusion coefficient is specified as $ D_o = 0.05 ,$ while in test 3, the dimensionless dispersion coefficient is also set as $ D_o = 0.05 .$ As described in test 1, after some trial and error, the backward time step is chosen as $ \Delta s' = 10^{-4} ;$ likewise, in test 3, $ \Delta t = 10^{-4} . $ 
Finally, in tests 1 and 2, the number of RW's launched from a chosen (forward-time) response point, $ \left( \mathbf{x}, s=0 \right) , $ is $ N_{\mathbf{x}, s=0} = 10^6 ;$ in test 3, $ N_{\mathbf{x'}, t'=0} = 10^6 , $ corresponding to the number of RW's launched from a (forward-time) impulse point, $ \left( \mathbf{x'}, t'=0 \right) . $ 

Based on these choices, and as detailed in Appendix A, the estimated relative spatial variation in captured RW's, $ \sigma_n / N , $ over the chosen forward solution interval, $ 0 < t \leq 10 , $ ranges between
\begin{equation*}
\sim 1.6 \times 10^{-5} < \frac{ \sigma_n}{N} < \hspace{0.2cm} \sim 1.6    
\end{equation*}
where the maximum variation occurs at late forward (early backward) time, $ t =\Delta t = 10^{-4} $ $ \left( s' = \Delta s' = 10^{-4} \right) . $ From Eq. (\ref{nrelvariation}) Appendix A, large, early-backward-time variations can be suppressed by choosing sufficiently large $ N = N_{\mathbf{x} , s=0} $ (or for forward time random walker integration, large enough $ N_{\mathbf{x'} , t'=0}).$ Our choice for $ N_{\mathbf{x} , s=0} $ $ \left( N_{\mathbf{x'} , t'=0} \right) $ reflects a compromise between computational cost and variance reduction. In addition, our chosen $ \Delta x' = \Delta y' = \Delta x = \Delta y = 10^{-2} $ are approximately five times the (average) incremental diffusion distance, $ \sqrt{D_o \Delta s'} \approx 2 \times 10^{-3} ; $ thus, and consistent with the results below, we appear to be obtaining statistically reasonable estimates of $ n \left( \mathbf{x'} , s' \right) $ $ [ n \left( \mathbf{x} , t \right) ] .$  

\subsection{Test Case 1, Part A: Deleterious depletion of random walkers at Dirichlet boundaries, remediation via random walker respawning, and error due to overlarge time step size}\label{respawnsection}
\begin{figure}[h]
\centering
\includegraphics[scale=0.45]{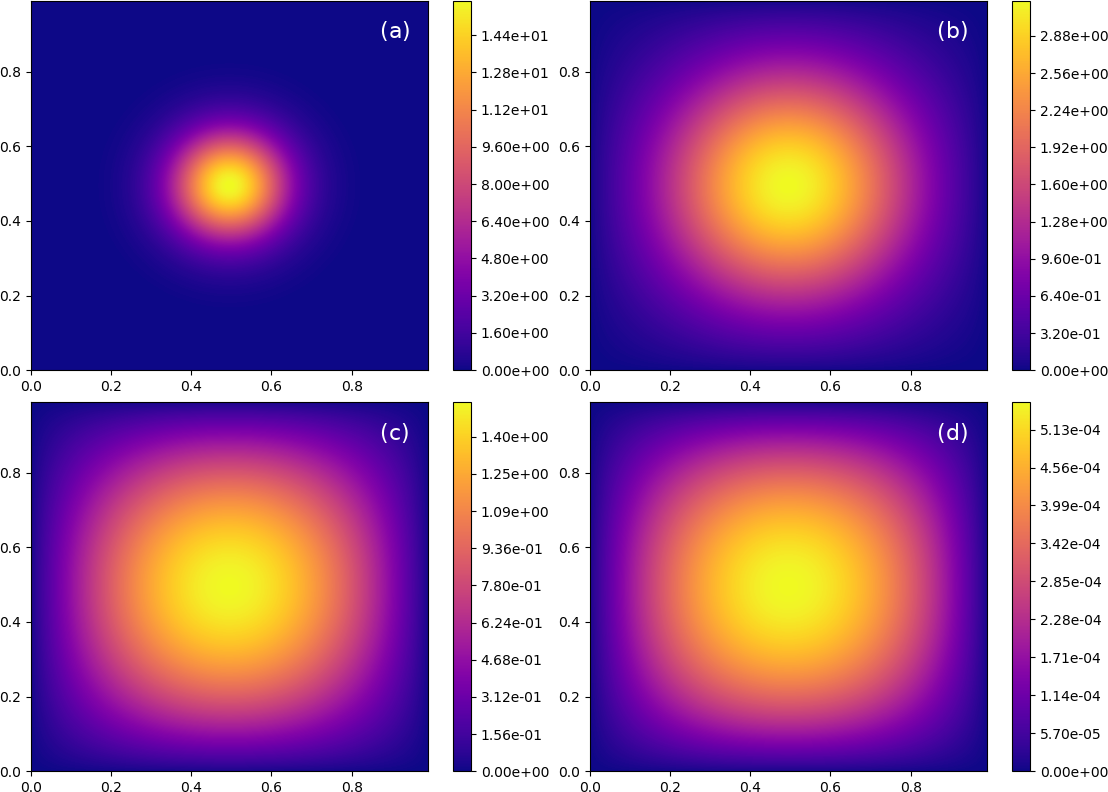}
\caption{Test 1: Exact Green's function, obtained using an analytical solution \cite{barton1989elements}, shown at forward time instants: (a) \(t' = 9.9,\) (b) \(t' = 9.5,\) (c) \(t' = 9.0,\) and (d) \(t' = 1.0.\) The total solution time is t=T=10.}
\label{fig:2}
\end{figure}

We use the first test case to emphasize the experimental nature of constructing Green's function estimators, briefly highlighting the challenges we encountered, as well as our solutions. In test 1, the solution (response) point is chosen as $ \left( \mathbf{x} , s \right) = \left( x=0.5, y=0.5, s=0 \right) . $ In all figures, we show estimated and actual GF's at various forward time instants, $ t' ,$ where $ 0 \leq t' \leq T=10 .$  For reference, Fig. \ref{fig:2} shows the exact Green's function.  

\begin{figure}[h]
\centering
\includegraphics[scale=0.45]{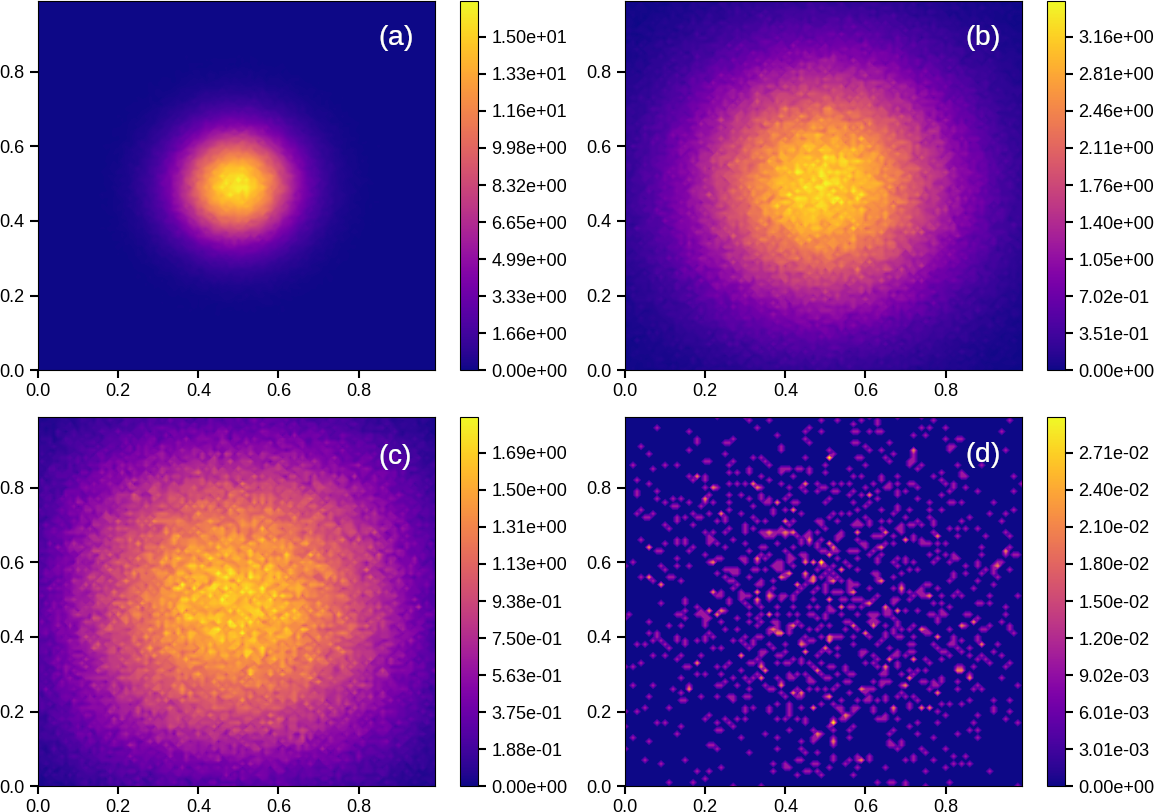}
\caption{Test 1, Part A: Backward time estimated Green's functions, shown at forward times (a) \(t' = 9.9,\) (b) \(t' = 9.5,\) (c) \(t' = 9.0,\) and (d) \(t' = 1.0.\) While the backward time-step, $ \Delta s' = 10^{-1} , $ is approximately 2 orders of magnitude larger than the recommended $ \Delta s_o' \approx \Delta A' / D_o = 2 \left( 10^{-3} \right) $, see Appendix A, the major source of error arises from  depletion of RW's at absorbing Dirichlet boundaries. Maximum relative errors in each estimate are: (a) $ e_{max} = 4.1  \% , $ (b) $ e_{max} = 10  \%, $ (c) $ e_{max} = 21  \%, $ and (d) $ e_{max} = 5200  \% .$ }
\label{fig:1}
\end{figure}

The detrimental effect of random walker absorption at Dirichlet boundaries is illustrated in Fig. \ref{fig:1}, where depletion of RW's becomes increasingly apparent as backward time integration of random walker trajectories proceeds. Comparing Figs. \ref{fig:2} and \ref{fig:1}, we observe that the maximum local error, 
\begin{equation}
e_{max} = max \  | G_{estimated} - G_{exact} |  / G_{exact} 
\end{equation}
increases from $ 4.1 \% $ at forward time, $  t' = 9.9 , $ to $ 5200 \% $ at $ t'=1 . $ Moreover, while  relatively late forward time estimates, obtained at $ t'=9.9, $ $ t'=9.5 ,$ and $ t' =9.0 , $ shown respectively in Figs. 2a, 2b, and 2c, remain qualitatively consistent with the actual Green's function in Fig. 1a-c, the early forward time estimate at $ t'=1 , $ Fig. 2d, is spatially amorphous, bearing no resemblance to the actual Green's function depicted in Fig. 1d. While the results in Fig. \ref{fig:1} were obtained using a time step, $ \Delta s' = 10^{-1} , $ which is approximately 2 orders of magnitude larger than the recommended $ \Delta s_o' \approx \Delta A' / D_o = 2 \left( 10^{-3} \right) $, see Appendix A, similar qualitative results, not shown, are observed for $ 10^{-4} \leq \Delta s \leq 10^{-1} . $


\begin{figure}[h]
\centering
\includegraphics[scale=0.45]{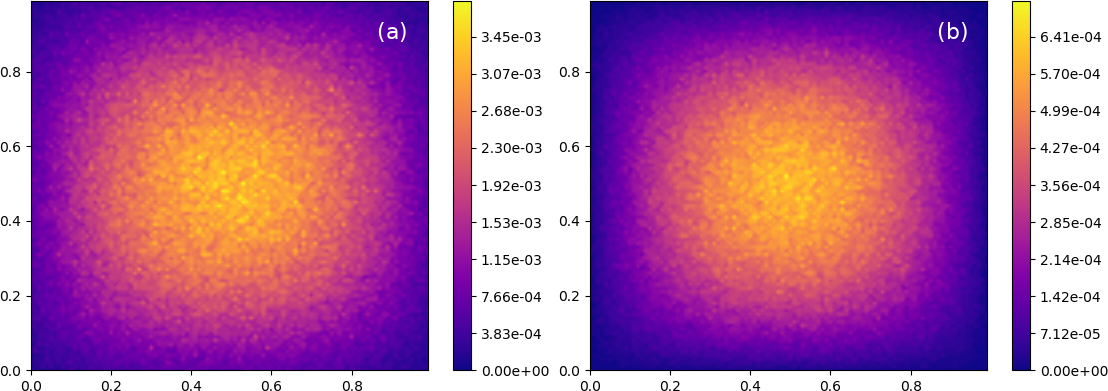}
\caption{Test 1, Part A: Estimated Green's function at \(t' = 1.0\), $ s' = 9 , $test 1, constructed using the respawning algorithm; compare with the exact solution in Fig. 1d and the estimate obtained in Fig. 2d, without respawning. The backward time step used in (a) is fifty times larger than the recommended time step, $ \Delta s_o' = \Delta A' / D_o , $ see Appendix A, and in (b), twenty times smaller than $ \Delta s_o' . $ Maximum errors, relative to the exact solution in Fig. 1d, are: (a) $ e_{max} = 570 \%, $ and (b) $ e_{max} = 25 \%.$ } 
\label{fig:3}
\end{figure}

Introduction of random walker respawning significantly improves both the accuracy and the qualitative consistency of Green's function estimates, as shown in Figs. \ref{fig:3}a and \ref{fig:3}b.  These depict, respectively, $ G_{est} (\mathbf{x},t=10|\mathbf{x'},t'=1) ,$ $ \mathbf{x} = \left( 0.5 , 0.5 \right) ,$ at the early forward time, $ t' = 1 $ (corresponding to $ s' =9 ).$ The two estimates shown in Figs. \ref{fig:3}a and \ref{fig:3}b use backward time steps, $ \Delta s' = 10^{-1} $ and $ \Delta s' = 10^{-4} , $ respectively, which bracket the recommended $ \Delta s' = \Delta s_o' \approx \Delta A' / D_o = 2 \left( 10^{-3} \right) .$ Again, similar qualitative results are observed for $ 10^{-4} \leq \Delta s \leq 10^{-1} . $


In order to improve the accuracy of the estimate shown in Fig. \ref{fig:3}a, the too-large time step, $ \Delta s' = 10^{-1} = 50 \times \Delta A' / D_o , $ is reduced to  $ 10^{-4} = 0.05 \times \Delta A' / D_o , $ reducing the maximum relative error, again observed at $ t'=1 ,$ from 570 \% to 25 \%.  
Although not shown, similar results are observed using $ \Delta s' = \Delta A' / D_o .$

\subsection{Test 1, Part B: Area averaging method for improving Green's function estimation accuracy and smoothness}\label{areavgsection}
While reducing $ \Delta s' $ to magnitudes on the order of $ \Delta s_{o}^{'} $ or less significantly improves Green's function estimation accuracy, relative error, at long backward integration times, here, $ s'= 9 , $ remains unacceptably large.  Moreover, as shown in Fig. \ref{fig:3}, plate b), long (backward) time estimated GF's exihibit spatial graininess, an artifact of using the naive estimator in Eq. (\ref{gengfstochasticconstructionrespawn}).   

In order to address both issues, we introduce a simple area-averaging technique in which the value of $ m \left( \mathbf{x}',s' \right) $ in (\ref{gengfstochasticconstructionrespawn}), obtained over an area, $ \Delta A' , $ encompassing grid point $ \mathbf{x}' = \left( x',y' \right), $ is replaced with the average of all $m$'s observed at all grid points lying within a square \((x\pm a, y\pm a)\), centered on $ \left( x',y' \right) .$ Here, the smoothing window size, \(a\), is the smaller of two values: the distance to the nearest boundary, or a maximum window size, $ a_{max} . $ At any given backward time instant, $ s' , $ $ a_{max} $ is determined by iteratively altering $ a $ in order to minimize the total root mean deviation, 
\begin{equation}\label{acalcn}
\overline{\sigma_G}= \frac{1}{M_{inc}} \sum_{i=1}^{M_{inc}'} \sqrt{ \left[ G(\mathbf{x}_i',s')  - \overline{G}(\mathbf{x}_i',s') \right]^2} 
\end{equation}
between individual estimated and averaged Greens functions. Here, $M_{inc}$ is the number of grid points included in the sum, $ G(\mathbf{x}_i',s') $ is the estimated Green's function at grid point $i$, as given by (\ref{gengfstochasticconstructionrespawn}), and $ \overline{G}(\mathbf{x}_i',s')$ is the averaged Green's function. As noted, we use known, exact GF's in place of $ \overline{G}(\mathbf{x}_i',s'). $ The prime on the sum specifies that only those grid points where the summand is not equal to zero are included. We include this additional filter since, in many cases, there are a significant number of points in both the analytical and numerical Green's function where both functions are zero. These flat areas, which can be seen in many late forward time GF's, can skew the total deviation to spuriously low magnitudes.

\begin{figure}[!h]
\centering
\includegraphics[scale=0.45]{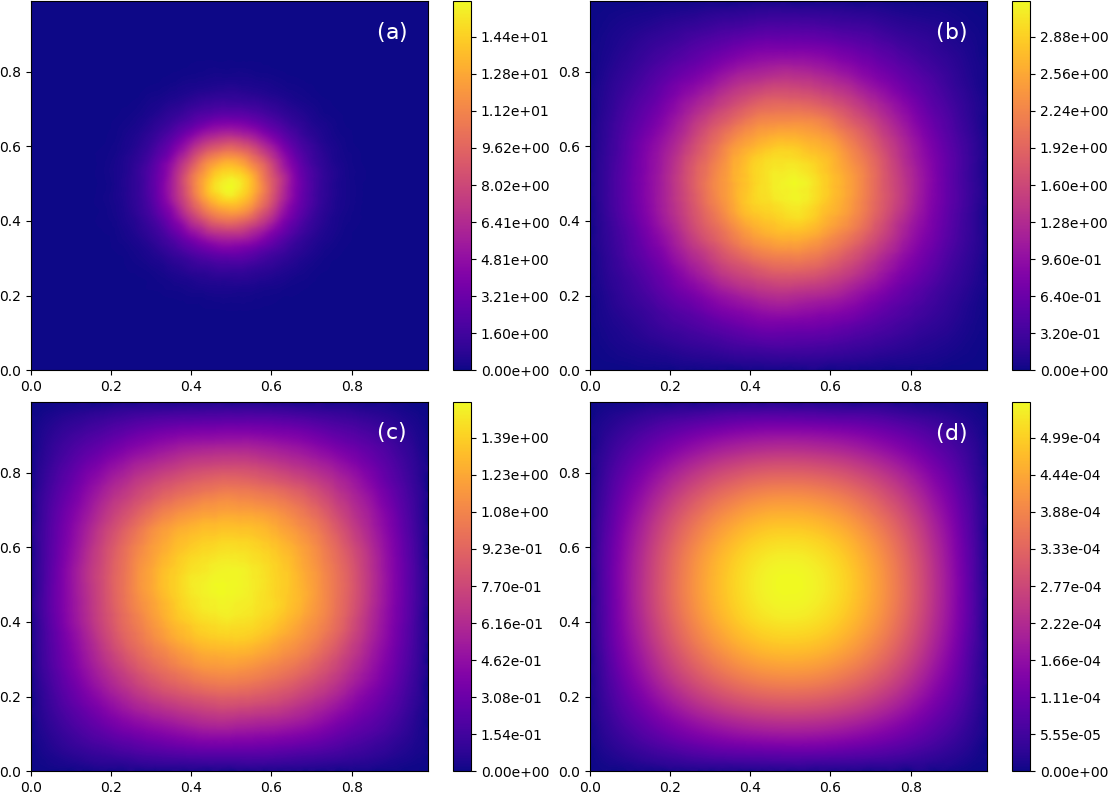}
\caption{Test 1, Part B: Estimated Green's functions, constructed using random walker respawning and area averaging. The backward time step, $ \Delta s' = 10^{-4} = 0.05 \times \Delta A' / D_o .$  Forward time instants, $ t' , $ maximum relative errors, $ e_{max} , $  and maximum smoothing window sizes, $ n_{max} = l'_s / \left( 2 \Delta x' \right) ,$ where $ l'_s $ is the dimensionless length of the square smoothing window, are, respectively: (a) \(t' = 9.9,\) $ e_{max} = 0.31 \% , $ $ n_{max} = 1,$ (b) \(t' = 9.5,\) $ e_{max} = 0.8 \%, $ $ n_{max} = 3,$ (c) \(t' = 9.0,\) $ e_{max} = 1.4 \% , $ $ n_{max} = 4,$ and (d) \(t' = 1.0,\) $  e_{max} = 2.9 \% , $ $ n_{max} = 13 . $ }
\label{fig:4}
\end{figure}

Finally, we observe that combining random walker respawning and area averaging provides relatively accurate estimated GF's, as shown in Fig.\ref{fig:4}. Comparing estimated and exact GF's, in Figs. \ref{fig:4} and \ref{fig:2}, respectively, it is clear that, even at large backward estimation times, $ s' = 9 , $ Fig. \ref{fig:4}d,  estimation accuracy and smoothness remain reasonable.  At small backward estimation times, $ s' = 0.1, $ $ s' =0.5 ,$ and $ s' = 1.0 ,$ shown respectively in Figs. \ref{fig:4}a - c, maximum local relative errors remain less than $ 1.5 \% ,$ and increase only to $ 2.9 \% $ at $ s'=9 ,$ in Fig. \ref{fig:4}d. Note too that Green's function magnitudes vary by more than five orders of magnitude over the chosen forward time interval, $ 1 \leq t' \leq 9.9 .$
Although not tested here, but as indicated by the $ \mathrm{O} \left( \Delta s' \right) $ weak convergence of the Milstein-Euler integration scheme, we believe that improved accuracy can be obtained, at significantly higher computational cost, by further reductions in $ \Delta s' . $  

Two final remarks are made concerning area averaging.

\vspace{0.2cm}

\noindent a) Beyond improved estimation accuracy and Green's function smoothing, area averaging reduces the number of random walkers required to achieve a given level of accuracy, significantly improving computational efficiency.  For example, experiments show that achieving $ \sim 1 \% $ late (backward) time accuracy, without area averaging, requires on the order of twenty million RW's.

\vspace{0.3cm} 

\noindent b) For poor Green's function estimates, like those shown in Figure \ref{fig:3}, plate a), no optimum window size can be found. There, the smoothing window minimizing the total deviation spans the entire domain, producing a constant function equal to the average of the values of all grid points. 

\vspace{0.4cm} 


\subsection{Test 2: Green's function estimates near curved boundaries}
\begin{figure}[h]
\centering
\includegraphics[scale=0.45]{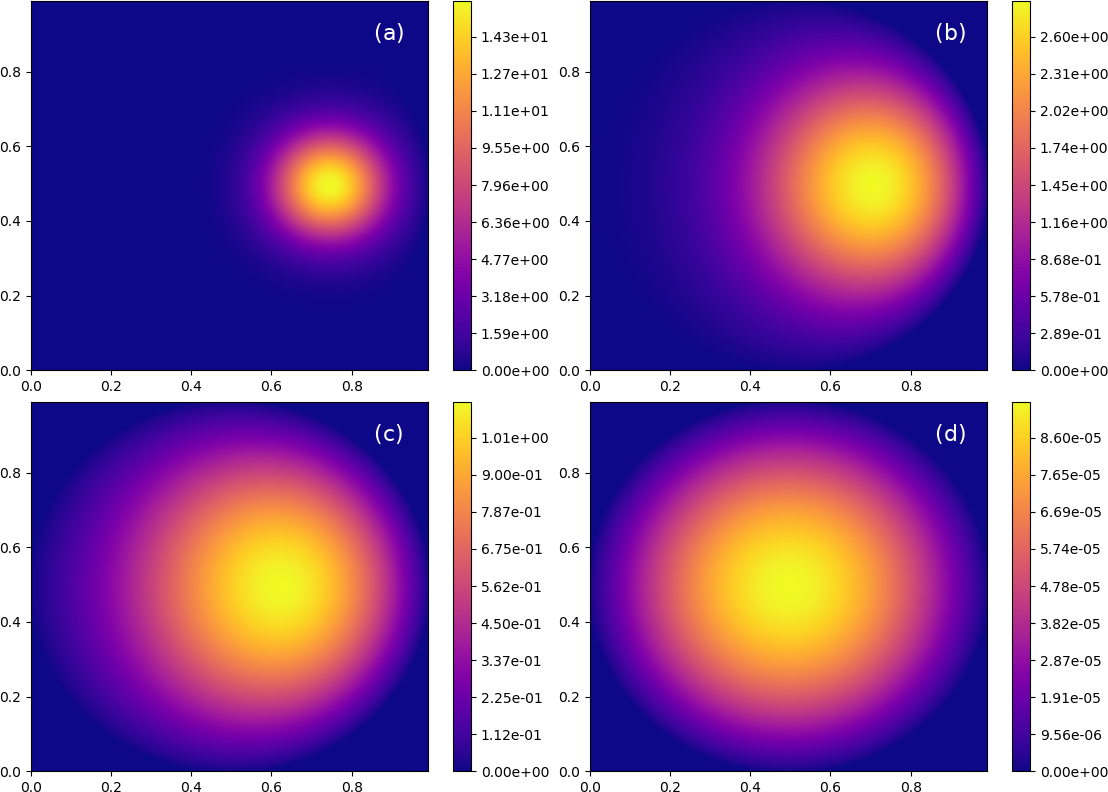}
\caption{Test 2: Exact Green's function in a circle of radius $ 0.5 , $ inscribed and centered within the unit square shown. The chosen solution (response) point is $ \left( x = 0.75, y = 0.5, t = 10 \right), $ and the exact Green's function is shown at forward time instants: (a) \(t' = 9.9,\) (b) \(t' = 9.5,\) (c) \(t' = 9.0,\) and (d) \(t' = 1.0.\)}
\label{fig:6}
\end{figure}

\begin{figure}[h]
\centering
\includegraphics[scale=0.45]{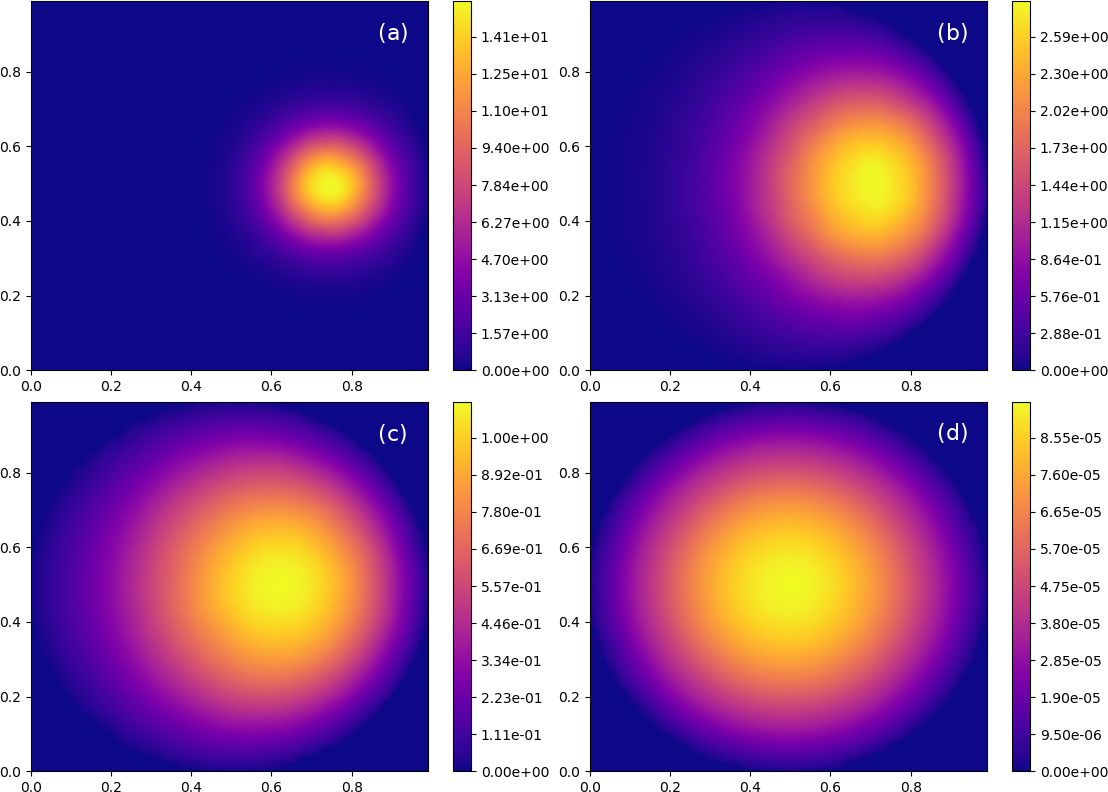}
\caption{Test 2: Estimated, area-averaged Green's functions in a circle of radius $ 0.5 , $ inscribed and centered within the unit square shown.  Random walker respawning is used and smoothing windows are circular. The unit square is discretized into $ 100 \times 100 $ random walker interrogation areas, $ \Delta A' = \Delta x' \times \Delta y' = 10^{-4} , $ introducing an $ \mathrm{O} \left( \Delta x' \right) $ discretization error in Green's function estimates. The solution point is at $ \left( x = 0.75, y = 0.5 \right). $ 
Forward time instants, $ t' , $ maximum relative errors, $ e_{max} , $  and maximum smoothing window sizes, $ n_{max} = a' / \left( \Delta x' \right) ,$ where $ a' $ is the radius of the circular smoothing window, are, respectively: 
(a) \(t' = 9.9,\) $ e_{max} = 1.57 \% , $ $ n_{max} = 2,$ 
(b) \(t' = 9.5,\) $ e_{max} = 0.86 \%, $ $ n_{max} = 4,$
(c) \(t' = 9.0,\) $ e_{max} = 0.89 \% , $ $ n_{max} = 6,$ 
and (d) \(t' = 1.0,\) $  e_{max} = 0.70 \% , $
$ n_{max} = 17.$}
\label{fig:5}
\end{figure}

Many evolution problems take place in regions having complicated, non-planar boundaries.  As a preliminary test of the stochastic Green's function estimation procedure in such problems, we perform a second set of numerical experiments in which GF's are constructed within a two-dimensional circular domain. Referring to Figures \ref{fig:5} and \ref{fig:6}, a circular domain (boundary not shown), having radius 0.5 and center at $ \left( x, y \right) = \left( 0.5 , 0.5 \right) $, is inscribed and centered within the unit square shown. Here, any random walker that reaches or passes through the circle, $ \left( x - 0.5 \right)^2 + \left( y - 0.5 \right)^2 = 0.25 , $ is absorbed, triggering random walker respawning at a random location within the circle. For simplicity, the square is discretized into the same uniform set of square $ 10^4 $ interrogation areas, $ \Delta A' \left( \mathbf{x'}, s' \right) = \Delta x' \Delta y' ,$ used in test 1, introducing an $ \mathrm{O} \left( \Delta x \right) $ discretization error in boundary-adjacent Green's function estimates. As in Example 1,  there is no flow \((\mathbf{b} = 0)\), transport takes place by isotropic diffusion, with $ D_o = 0.05; $ we use a backward-time step, $ \Delta s' = 10^{-4} = 0.05 \times \Delta A' / D_o ,$ and as in test 1, $ N_{\mathbf{x},s } = 10^6 $ RW's are launched, in this case, from a right-of-center solution point $ \left( \mathbf{x}, t \right) = \left( x=0.75, y=0.5, t=T=10 \right) .$ 

In order to adapt the smoothing procedure to a circular domain, we define the smoothing window as a circle of radius \(a\) centered on the grid point of interest. As before, \(a\) is allowed to vary, being the smaller of two values: the distance to the boundary and the maximum window size. By this approach, the smoothing window is always contained within the circular domain. Although not explored here, beyond circular and square/rectangular smoothing window shapes - the former useful for domains having curved boundaries - polygonal shapes may also prove useful.

Exact and estimated GF's are shown respectively in Figs. \ref{fig:6} and \ref{fig:5}. In this case, and in contrast to test 1, the largest estimation errors are observed at the earliest backward time instant, $ \Delta s' = 0.1 . $ Comparing patterns in maximum relative error evolution in test 1 and 2, Figs. \ref{fig:4} and \ref{fig:5}, as expected, in test 1, $ e_{max} $  progressively increases with increasing $ s' , $ at least for $ s' = 0.1, $ $ s'= 1.0, $ and $ s' = 9.0 .$ By contrast, in test 2, $ e_{max} $ progressively decreases at all four probed time instants. These results indicate the following: I) The fact that $ e_{max} $'s at $ s'=1.0 $ and $ s'=9 ,$ in test 1, Fig. \ref{fig:4}, plates c) and d), are significantly larger than those in test 2, Fig. \ref{fig:6}, plates c) and d), suggests that use of circular smoothing windows in test 2, as opposed to square windows in test 1, reduces estimation error, even under conditions where near-boundary Green's function estimates are subject to an $ \mathrm{O} \left( \Delta x' \right) $ discretization error. II) In light of point a), the relatively large error at $ s' =0.1 $ in test 2 appears to reflect our choice of the asymmetrically placed solution point, $ (x=0.75, y=0.5).$ 

Crucially, over the forward time span, $ 0.1 \leq t' \leq 9 , $ high-fidelity Green's function estimates are obtained, here, over approximately six orders of magnitude in local Green's function magnitudes.

\section{Test 3: Forward time Green's function estimation - advection-dispersion-reaction problems}\label{test3section}
The third example addresses two questions:  

\vspace{0.2cm}

\noindent a) Can Green's functions be stochastically estimated using forward time random walker integration? 

\vspace{0.2cm}

\noindent b) How well does the stochastic Green's function estimation procedure perform for advection-dispersion (or diffusion)-reaction problems?

As noted, Green's function estimation using forward time random walker integration corresponds to estimating the \textit{response} at a set of $ Q $ space-time points, $ \left[ \mathbf{x_1} , t; \mathbf{x_2} , t; ... \mathbf{x_Q} , t, \right] ,$ produced by a (delta function) impulse at a single point, $ \left( \mathbf{{x}'} , t' \right) : $ $ G \left( \mathbf{x}, t | \mathbf{{x}'} , t' \right) .$ By contrast, Green's function estimated via backward time integration provides the response at the point, $ \left( \mathbf{x}, t \right) , $ as produced by a (discrete) set of $ Q $ impulse points, $ \left[ \mathbf{x_1'} , t'; \mathbf{x_2'} , t'; ... \mathbf{x_Q'} , t' \right].$ Thus, for physical transport problems where property $ \eta $ is not required over an entire problem domain, backward integration offers distinct cost advantages. 

In order to obtain an analytical Green's function, we take advantage of an exact solution obtained by \cite{sanskrityayn2018analytical} for a ground water advection-dispersion-reaction problem, governed by 
\begin{equation}\label{groundwater}
\nabla^T \otimes \nabla : \mathbf{D} \hspace{0.05cm} \eta - \nabla \cdot \left( \mathbf{v} \eta \right) - \gamma\eta - \frac{\partial\eta}{\partial t} = -f(\mathbf{x}, t),
\end{equation}
corresponding to a slightly generalized version of Eq. (\ref{mainequation1}), and designed to capture a wide range of ground water transport conditions \cite{sanskrityayn2018analytical}, including, e.g., variable density flow in near-coastal aquifers \cite{bearcheng}. 

Importantly, using the substitutions preceding Eq. (\ref{mainequation}), the sink term, $ \gamma \eta ,$ in Eq. (\ref{groundwater}) can be suppressed, so that the transformed version of Eq. (\ref{groundwater}) has the generic form of a Fokker-Planck (FP) equation, with $ p \left( \mathbf{x}, t | \mathbf{x'}, t' \right) $ replacing $ \eta \left(\mathbf{x}, t \right) :$
\begin{equation}\label{fptest3}
\nabla^T \otimes \nabla : \mathbf{D} \hspace{0.05cm} p - \nabla \cdot \left( \mathbf{v} p \right) -  \frac{\partial p}{\partial t} = 0,
\end{equation}
and where the delta function source term treated in \cite{sanskrityayn2018analytical}, 
\begin{equation}
    f \left( \mathbf{x}, t \right) = \delta \left( \mathbf{x} - \mathbf{x'} \right) \delta \left( t - t' \right) 
\end{equation}
corresponds to the same-time condition,
\begin{equation}\label{sametimeFP}
p \left( \mathbf{x} , t | \mathbf{x'},t' \right) \rightarrow \delta \left( \mathbf{x} - \mathbf{x'} \right) \hspace{0.7cm} t \rightarrow t'
\end{equation}
imposed on the FP equation.  Thus, the analytical solution for $ \eta \left(\mathbf{x}, t \right) $ obtained in \cite{sanskrityayn2018analytical} can be interpreted as a solution to the Fokker-Planck equation for $p \left( \mathbf{x} , t | \mathbf{x'},t' \right).$

We show in Appendix D that the transport equation solved by the estimated forward time Green's function is given by
\begin{equation}\label{mainequationfwd}
\mathbf{D} \mathbf{:} \nabla^T \otimes \nabla \eta - \mathbf{v_a}\cdot\nabla\eta - \frac{\partial\eta}{\partial t} = -f(\mathbf{x}, t),
\end{equation}
where the actual advective flow field, $ \mathbf{v_a} \left( \mathbf{x} , t \right) ,$ in (\ref{mainequationfwd}) corresponds to a reversed flow field in Eq. (\ref{fptest3}): $ \mathbf{v} \left( \mathbf{x} , t \right) = - \mathbf{v_a} \left( \mathbf{x} , t \right) .$ 
\subsection{Recipe for forward time Green's function estimation, in brief}
Since the Fokker-Planck equation governs the transition density for observing members of a random walker swarm at response point $ \left( \mathbf{x}, t \right),$ launched in the forward time direction from impulse point $ \left( \mathbf{x'}, t' \right),$ computationally, we can use the same recipe given above for estimating $ p \left( \mathbf{x} , t | \mathbf{x'},t' \right) $ via backward random walker integration. Here, we highlight necessary modifications. 

The stochastic differential equation governing the forward time motion of individual RW's is given by 
\begin{equation}\label{stocheqnfwd1}
d \boldsymbol{\chi} \left( t \right) = \mathbf{v} \left( \boldsymbol{\chi} \left( t \right) ,t \right) dt + \mathbf{B} \left( \boldsymbol{\chi} \left( t \right),t \right) \cdot d \mathbf{w} \left( t \right)
\end{equation}
where, again, $ \mathbf{w} \left( t \right) $ is an n-dimensional Weiner process, and $ \mathbf{v} \left( \mathbf{x} , t \right) $ is the velocity field. 

In forward time integration, a swarm of $ N_{\mathbf{x'},t'} $ is launched from impulse point $ \left( \mathbf{x'}, t' \right) ,$ where $  t \geq t' ,$ with integration stopped at a chosen final time, $ t =T .$ The spatial solution domain is discretized into $ N_{\Delta A}  $ interrogation areas, $ \Delta A \left( \mathbf{x} , t \right) = \Delta x \Delta y ,$ and the number of RW's reaching each $ \Delta A \left( \mathbf{x} , t \right) , $ $ n \left( \mathbf{x}, t \right) ,$ at a set of chosen times, $ t_1 , t_2 , ..., t_q = T , $ is tabulated.   
As noted immediately below, random walker respawning is not needed in Test 3. Thus, estimated GF's follow from Eq. (\ref{stochestgf}):
\begin{equation}\label{stochestgfFP}
G \left( \mathbf{x} , s | \mathbf{x'},s' \right) = p \left( \mathbf{x} , s | \mathbf{x'},s' \right) = \frac{n \left( \mathbf{x}, t \right)}{N_{\mathbf{x'}t'} \Delta A \left( \mathbf{x}, t \right)} +  O \left( {\Delta t} \right)
\end{equation}

Following \cite{sanskrityayn2018analytical}, the components of the dispersion tensor, $ \mathbf{D} , $ are given by $ D_{xx} = D_0 \left[a_1 x + a_2)^2\psi_1 \left(mt\right) \right], $ $ D_{yy} = D_0 \left[b_1 y + b_2)^2\psi_2 \left(mt \right)\right] , $ with off-diagonal terms, $ D_{xy} = D_{yx} = 0 . $ Note again that stochastic estimation of Green's function is limited to diffusion and dispersion problems characterized by diffusive or dispersive transport dominant in one, two, or three mutually orthogonal directions. Following \cite{sanskrityayn2018analytical}, the ground water flow field is given by $ v_x = v_0 \left(a_1 x + a_2)\psi_1 \left(mt \right)\right) $ and $ v_y = v_0 \left(b_1 y + b_2\right)\psi_2 \left(mt \right),$ with $ \psi_1 \left(mt \right) $ and $ \psi_2 \left(mt \right) $ being two arbitrary functions of time. For purposes of this test, we set \(D_0 = 0.05\), \(a_1 = b_1 = 1\), \(a_2 = b_2 = 0\), \(v_0 = 0.2\) and \(\psi_1 \left(mt \right) = \psi_2 \left(mt\right) = 1\). Finally, in  the analytical solution presented by \cite{sanskrityayn2018analytical}, we set the depletion constant, $ \gamma = 0.5. $

\begin{figure}[h]
\centering
\includegraphics[scale=0.45]{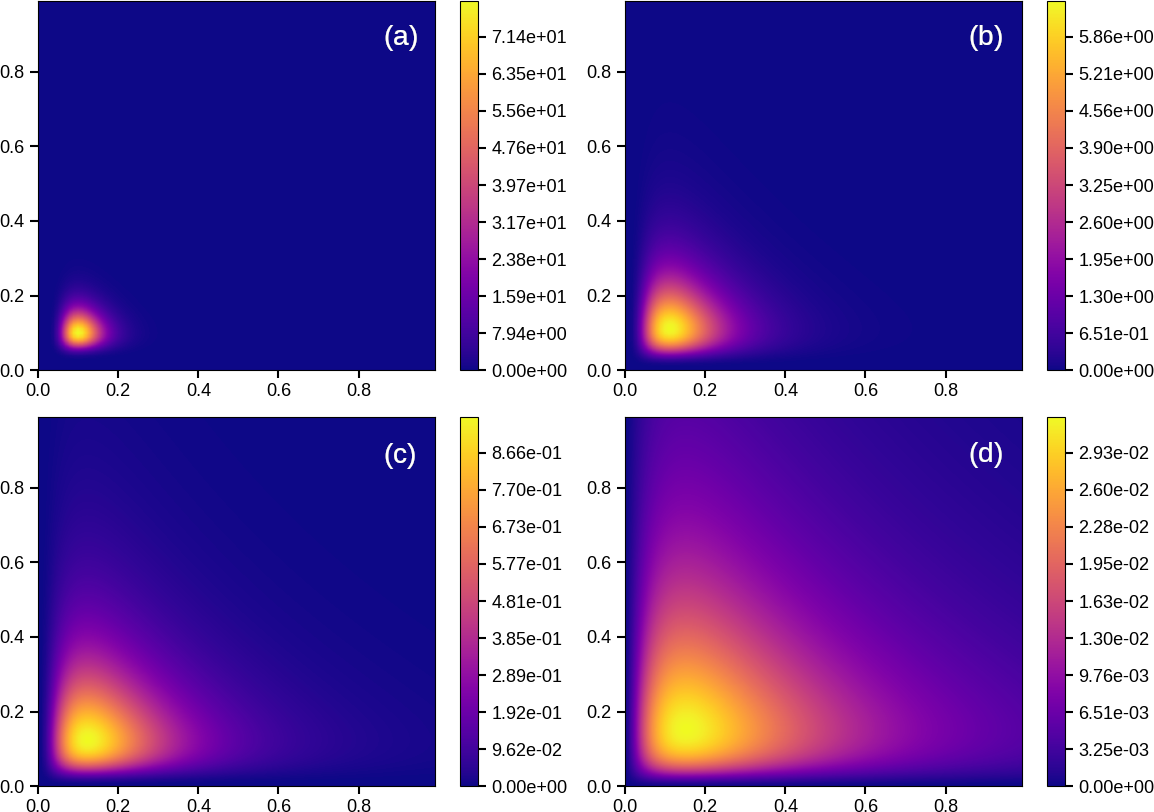}
\caption{Test 3: Exact Green's functions on a semi-infinite domain, applicable to advective-dispersive-reactive ground water transport in a semi-infinite, anisotropic region. Here, the advective flow field is directed at $ 45^o $ from horizontal and accelerates linearly at a fixed rate, both from left to right and from bottom to top. Green's functions, constructed using an exact ground water transport solution \cite{sanskrityayn2018analytical}, are shown at: (a) \(t' = 1,\) (b) \(t' = 3,\) (c) \(t' = 5,\) and (d) \(t' = 9.\)}
\label{fig:aflow}
\end{figure}

\begin{figure}[h]
\centering
\includegraphics[scale=0.45]{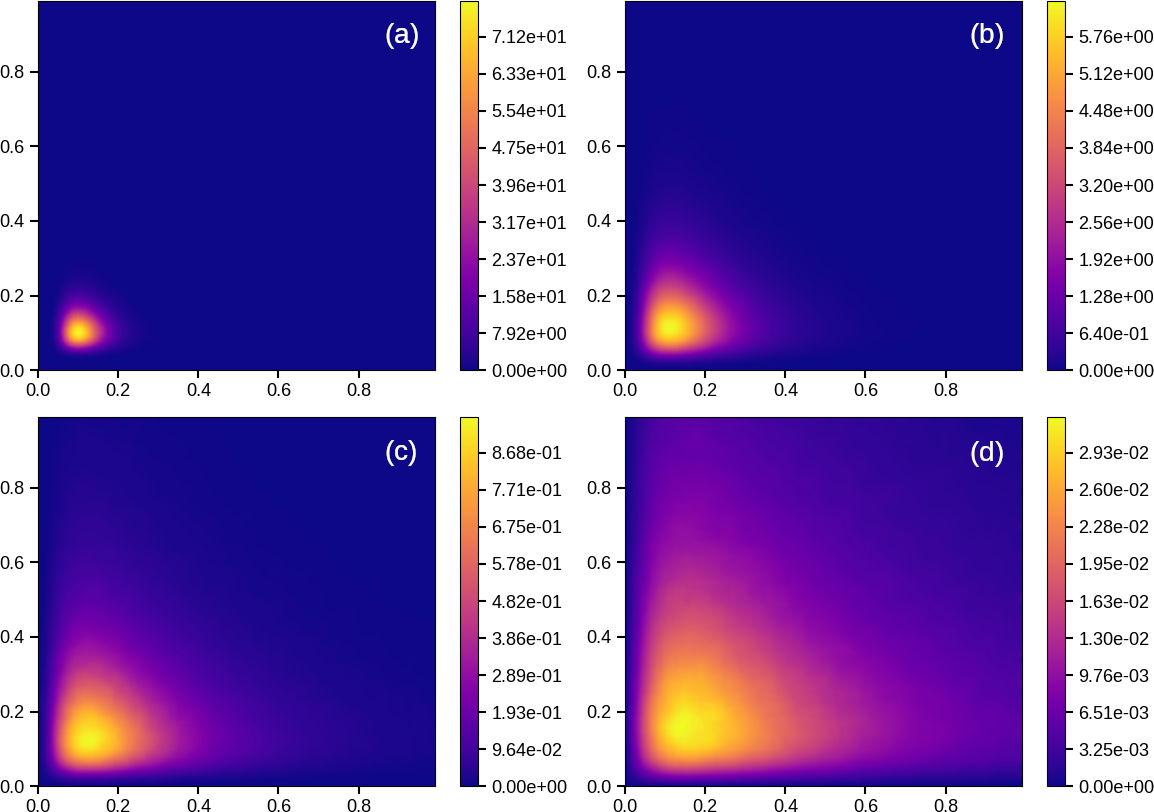}
\caption{Test 3: Estimated, smoothed Green's functions, constructed using a backward time step $ \Delta s' = 10^{-4} .$ Forward time instants, $ t' , $ maximum relative errors, $ e_{max} , $  and maximum smoothing window sizes, $ n_{max} = l_s' / \left( 2 \Delta x' \right) ,$ where $ l_s' $ is the length of the square smoothing window, are, respectively:
(a) \(t' = 1.0,\) $ e_{max} = 0.63 \% , $ no smoothing, (b) \(t' = 3.0,\) $ e_{max} = 1.73 \%, $ $ n_{max} = 1,$ (c) \(t' = 5.0,\) $ e_{max} = 0.52 \% , $ $ n_{max} = 1,$ and (d) \(t' = 9.0,\) $  e_{max} = 0.15 \% , $ $ n_{max} = 2.$}
\label{fig:nflow}
\end{figure}


In test 3, two features allow us to bypass the random walker respawning algorithm, and in addition, reduce the solution domain.  First, since the physical problem domain is infinite, RW's, in reality, never reach a domain boundary.  In order to avoid specification of domain dimensions that are too small, allowing unphysical random walker absorption, various approaches can be used.  For example, one can choose boundary distances, $ d_{boundary} ,$ that are some greater-than-unity multiple of the characteristic advective and diffusive transport distance, $ v_{char} T + \sqrt{ D_{char} T } , $ where $ v_{char} $ and $ D_{char} $ are characteristic velocity and diffusion/dispersivity magnitudes, and $ T $ is the chosen total solution time. As described immediately below, we use a second approach, representing a refined, and likely more economical variation on this strategy. 

Second, since the ground water velocity field and the position-dependent dispersivity both tend to zero as the origin,  $ \left( x,y \right) = \left( 0,0 \right) $ is approached, very few RW's exit the first/positive quadrant of the infinite space. Thus, for this test, calculations are limited to the first quadrant and, based on numerical experiments showing that random walker respawning does not significantly improve Green's function estimation accuracy, respawning is turned off.

In order to define the  solution domain,  at the last forward time instant, $  t=T=10, $ we determine the average \(x-\) and \(y-\) positions of the random walker swarm launched from the origin, $ \mu_x \left( t=T \right) $ and $ \mu_y \left( t=T \right) ,$ as well as standard deviations in these averages, $ \sigma_x \left( t=T \right) $ and $ \sigma_y \left( t=T \right) . $ The non-normalized (but dimensionless) extent of the problem domain is then defined as: \(0 < x < \mu_x + \sigma_x\) and \(0 < y < \mu_y + \sigma_y\). The solution domain is discretized into $ 10^4 $ interrogation areas, but with different fixed sizes in the $x-$ and $y-$directions (since $\mu_x + \sigma_x \neq \mu_y + \sigma_y)$. As in tests 1 and 2, we choose $ \Delta t =  10^{-4} $ and $ N_{\mathbf{x'},t'} = 10^6 . $ Finally, the same smoothing technique used in tests 1 and 2 is applied here.

Stochastically estimated GF's are plotted in Figure \ref{fig:nflow} 
at the four forward time instants: $ t = 1,$ $ t=3 ,$ $ t=5 , $ and $ t=9 . $ For comparison, exact GF's, determined from \cite{sanskrityayn2018analytical}, are shown in Fig. \ref{fig:aflow}. Again, quality Green's function estimates are obtained at all four probed time instants, validating the use of forward time estimation.  Again, under most circumstances, however, backward time estimation is recommended since  one Green's function estimate, obtained at a single response point, $ \left( \mathbf{x}, t \right) , $ can be used in the magic rule to estimate the evolution of $ \eta \left( \mathbf{x}, t \right) , $ subject to any set of forcing conditions. 

\section{Discussion}
When they can be found, Green's functions provide a powerful tool for solving linear - and incrementally, nonlinear - evolution problems having the generic forms given by (\ref{mainequation1}) and (\ref{mainequation}). Unfortunately, finding GF's 
can be extremely difficult, particularly within geometrically complicated regions. 

This paper proposes a numerical method for stochastically estimating GF's, appropriate to diffusive and advective-dispersive/diffusive-reactive  transport of scalar- and vector-valued quantities, e.g., convective heat and mass transfer problems involving passive, as well as chemically-reacting species. We validate the technique against three known GF's, two that apply to simple diffusion in finite regions, and a third applicable to ground water transport in an infinite region. As formulated here, the stochastic Green's function estimation technique can be applied to scalar transport in independently determined turbulent flows.  More fundamentally, the technique might also be adapted to generate, in time-incremental fashion, GF's for computing turbulent flows. In this case: i) the approach is limited to conditions where fluid particle momentum and vorticity dynamics are strongly Brownian, ii) the transported quantity corresponds to either local fluid momentum or vorticity \cite{keanini2011green}, and iii) the turbulent flow problem is parabolic, evolving from a known initial condition. 

We identify three principle challenges and propose straightforward solutions for each:

\vspace{0.2cm}

\noindent a) Since the probability of observing RW's at Dirichlet boundaries is zero, RW's that reach such boundaries must be removed. Random walker absorption, in turn, quickly compromises the quality of estimated GF's.  In order to address this problem, we introduce a random walker respawning technique which conserves the number of RW's initially launched, as well as enforcing continuum mass conservation and diffusion at each respawning location.  

\vspace{0.2cm}

\noindent b) Estimated Green's functions exhibit spatial graininess, an artifact associated with use of discrete naive transition density estimators \cite{silverman}, Eq. (\ref{naiveestp}). In order to solve this problem, as well as improve Green's function estimation accuracy, we propose an area-averaging technique that bins neighboring Green's function estimates, taken over small smoothing windows encompassing any given interrogation area, $ \Delta A' .$ For most points, $ \mathbf{x'} , $ not near a boundary, the size of the smoothing window, $ a , $ is determined in iterative fashion by minimizing the point-wise deviation between numerical GF's and averaged GF's. In this first exploratory study, we use known rather than averaged GF's in determining $ a .$ 

\vspace{0.2cm}

\noindent c) To stochastically estimate Green's functions, as given by Eq. (\ref{stochestgf}) (no respawning, test 3) or (\ref{gengfstochasticconstructionrespawn}) (with respawning, tests 1 and 2), three numerical parameters must be chosen: i) the number of RW's, $ N_{\mathbf{x}, s} $ $ ( N_{\mathbf{x'}, t'} $ for forward time random walker integration), launched from the chosen solution/response point, $ \left( \mathbf{x}, s \right) $ [chosen impulse point, $\left( \mathbf{x'}, t' \right) ]$  ii) the backward time-step size, $ \Delta s' $ [forward time-step size, $ \Delta t ],$ and iii) the random walker interrogation area, $ \Delta A' $ $ [ \Delta A ].$  Appendix A outlines approaches for determining these parameters; to the best of our knowledge, and beyond \cite{schuss2}, the literature does not offer much guidance on this important practical question.   


\vspace{0.2cm}

Although not treated in this paper, the proposed stochastic Green's function construction technique provides a dual framework for tackling coupled particle and continuum transport problems, important, for example, in microscale bioheat and biomass transport problems  \cite{schuss2}, and a feature absent in BEM, GEM, and FEM formulations. The Green's function estimation recipe shows that this particle-continuum duality arises from the connection between stochastic differential equations modeling particle-like transport and associated  Fokker-Planck and backward Kolmogorov equations, probabilistically describing the same transport. 

\section*{Acknowledgements}
This work was supported by the Office of Naval Research, grant no. N00014-18-1-2754.  Technical and coding support provided by Nikita Nikulsen is gratefully acknowledged.

\bibliography{sampleold2}

\newpage

\section*{Supporting information}
\paragraph*{\textbf{S1}}
\label{S1}
{\bf Appendices to 'Stochastic estimation of Green’s functions with application to advection-diffusion-reaction problems'}



\vspace{0.3cm}

\paragraph*{\textbf{S1 Appendix A}}
\label{S1 Appendix A}
{\bf Choosing the number of RW's launched, $ N_{\mathbf{x'}, t'} ,$ interrogation area, $ \Delta A \left( \mathbf{x} , t \right) , $ and the time step, $ \Delta t. $} The stochastic estimation technique requires specification of three numerical parameters: i) the backward (forward) time-step size, $ \Delta s' $ $ \left( \Delta t \right) , $  ii) the interrogation area, $ \Delta A' \left( \mathbf{x}',s' \right) $ $ \left( \Delta A' \left( \mathbf{x}',s' \right) \right)  ,$ and iii) the number of RW's, $ N_{\mathbf{x},s} $ $ \left( N_{\mathbf{x'},t'}\right) , $ launched from the chosen GF solution point, $ \left( \mathbf{x}, s \right) $ $ \left( \left( \mathbf{x'}, t' \right) \right) .$ Since the literature appears to offer little guidance on choosing these parameters, and since the choices require some experimentation, we propose two straightforward approaches for identifying appropriate values. Here, we outline the technique used in this study. We focus on the simplest case in which the diffusion tensor is uniform and isotropic, $ \mathbf{D}= D_o \mathbf{I} . $ It is important to recognize that in both approaches, some numerical experimentation is required.


In order to address these important questions, we focus on forward-time integration of RW trajectories:
\begin{equation}\label{stocheqnfwd}
d \boldsymbol{\chi} \left( t \right) = \mathbf{v} \left( \boldsymbol{\chi} \left( t \right) ,t \right) dt + \mathbf{B} \left( \boldsymbol{\chi} \left( t \right),t \right) \cdot d \mathbf{w} \left( t \right)
\end{equation}
and highlight two features: a) the connection between $ \Delta A \left( \mathbf{x} , t \right) $ and $ \Delta t ,$ and b) minimization of early-time variance in the number of RW's captured on $ \Delta A \left( \mathbf{x} , t \right) .$ 

For simplicity, limit attention to two-dimensional diffusive transport of $ \eta \left( \mathbf{x}, t \right) .$ Spatially resolved random walker displacements, $ \mathbf{d} \boldsymbol{\chi} \left( t \right) , $ over $ \Delta t , $ require $ \Delta A \left( \mathbf{x} , t \right) \sim \langle d \mathbf{\chi} \left( t \right) \cdot d \mathbf{\chi} \left( t \right) \rangle = D_o \Delta t .$ Thus, we choose
\begin{equation}\label{areaest}
\Delta A \left( \mathbf{x} , t \right) = D_o \Delta t 
\end{equation}
where this choice represents the smallest area allowing statistically acceptable sampling of $ n \left( \mathbf{x},t \right) : $ If $ \Delta A \left( \mathbf{x}, t \right) $ is chosen smaller than $ D_o \Delta t , $ then for any given number of RW's launched from $ \left( \mathbf{x'}, t' \right) ,$ $ N_{\mathbf{x'}, t'} , $ a number of RW's located near $ \Delta A \left( \mathbf{x} , t \right) ,$ at forward instant $ t - \Delta t , $ will (numerically) jump over/around $ \Delta A \left( \mathbf{x} , t \right) $ at $ t .$ The result is an underestimated $ n \left( \mathbf{x},t \right) . $  Likewise, if $ \Delta A \left( \mathbf{x},t \right) $ is chosen larger than $ D_o \Delta t , $ the number captured on $ \Delta A \left( \mathbf{x} , t \right) $ is too large, overestimating $ n \left( \mathbf{x},t \right) . $ 

In order to determine acceptable $ N_{\mathbf{x'},t'} ,$ consider the time-dependent standard deviation, $ \sigma_n , $ in $ n \left( \mathbf{x},t \right) , $ in the number of RW's captured by the set of interrogation areas, $ \{ \Delta A_1 \left( r \left( \tau \right) \right) ,  \Delta A_2 \left( r \left( \tau \right) \right) , ... \Delta A_{N_{\Delta A}} \left( r \left( \tau \right) \right)  \} , $ lying on the circle, $ r = \sqrt{ D_o \tau} ,$ centered on the launch location, $ \mathbf{x'} : $
\begin{equation}\label{stddevn}
\sigma_n \left( \tau \right) = \sqrt{\frac{1}{N_{\Delta A} \left( r \right)} \sum_{i=1}^{N_{\Delta A} \left(r\right)} \left[ n_i \left( r \right) - \overline{n} \left( r \right) \right]^2 }
\end{equation}
where $ r = r \left( \tau \right) $ is the mean penetration radius of the $ N_{\mathbf{x'} , t ' } $ RW's launched from $ \mathbf{x'} $ at time $ t' , $ $ \tau = t - t' ,$ and $ \overline{n} \left( r \right) = N_{\Delta A}^{-1} \sum_{i=1}^{N_{\Delta A}} n_i \left( r \right) .$

The number of interrogation areas lying on $ r \left( \tau \right) ,$ $ N_{\Delta A} ,$ is given approximately by:
\begin{equation}
N_{\Delta A} \sqrt{ \Delta A_{\mathbf{x}, t }} \sim r \left( \tau \right) 
\end{equation}
where, in the case that the the spatial domain is discretized into uniform interrogation areas, $ \Delta A_{\mathbf{x}, t } = \Delta x \cdot \Delta y ,$ with $ \Delta x \approx \Delta y ,$ $ \sqrt{ \Delta A_{\mathbf{x}, t }} \sim \Delta x \approx \Delta y . $ Since $ \Delta x \approx \Delta y \sim \sqrt{ D_o \Delta t } ,$ then
\begin{equation}
N_{\Delta A} \sim \frac{ \sqrt{D_o \tau }}{\Delta x } \sim \frac{\sqrt{D_o \tau}}{\sqrt{D_o \Delta t}} = \sqrt{\frac{\tau}{\Delta t}}
\end{equation}

From Eq. (\ref{stddevn}), and for small interrogation areas, $ \Delta A_{\mathbf{x},t} << A_o ,$ $A_o$ being the area of the solution domain,
\begin{equation}\label{nstddev}
\sigma_n \left( \tau \right) \sim \frac{1}{ \sqrt{ N_{\Delta A}} } \sim \left( \frac{ \Delta t}{\tau} \right)^{1/4}
\end{equation}
From this expression, it becomes apparent that the maximum spatial variation, $ \sigma_n \left( \tau \right),$ in random walker capture numbers, $ n \left( \mathbf{x}, t \right) ,$ occurs at $ \tau_{min} = \Delta t ,$ where $ \sigma_n \left( \Delta t \right) \sim 1 .$

In order determine an acceptable $ N_{\mathbf{x'},t'} ,$ we calculate the ratio of the mean number of RW's captured at $ \tau = \Delta t ,$ $ \langle n \left( \mathbf{x}, \Delta t \right) \rangle , $ to $ N_{\mathbf{x'} , t' } :$
\begin{align}\label{fractionntoN}
\frac{\langle n \left( \mathbf{x}, \Delta t \right) \rangle }{N_{\mathbf{x'} , t' } }= P \big[ \boldsymbol{\chi} \left( \Delta t  \right) \in \Delta A \left( \mathbf{x}, \Delta t \right) | \boldsymbol{\chi} \left( t' \right) = \mathbf{x'} , t' \big] \nonumber \\
 &  \hspace{-7cm} = \int_{\Delta A \left( \mathbf{x}, \Delta t \right) } \frac{1}{4 \pi D_o \Delta t } \exp \bigg[ \frac{ - \left( \mathbf{x} - \mathbf{x'} \right)^2}{4 D_o \Delta t } \bigg] d \mathbf{x}
\end{align}
or 
\begin{equation}\label{fractionntoNfinal}
\frac{\langle n \left( \mathbf{x_j}, \Delta t \right) \rangle }{N_{\mathbf{x'} , t' } }= \frac{1}{ \pi } \exp \bigg[ \frac{ - \left( \mathbf{x_j} - \mathbf{x'} \right)^2}{4 D_o \Delta t } \bigg] \frac{\Delta x \Delta y}{4 D_o \Delta t} \sim \frac{ \Delta x \Delta y }{ 4 \pi e^{1/4} D_o \Delta t}
\end{equation}
where $ \mathbf{x_j} $ represents, e.g., the centroid of the $ j^{th} $ interrogation area, $ \Delta A_{\mathbf{x_j} , \Delta t} ,$ lying on $ r \left( \Delta t \right) = D_o \Delta t ,$ and where $ \left( \mathbf{x_j} - \mathbf{x'} \right)^2 \sim D_o \Delta t .$

Finally, solving for $ \langle n \left( \mathbf{x_j}, \Delta t \right) \rangle , $ we obtain the maximum relative spatial variation in captured RW's, expressed as a function of $ N_{\mathbf{x'},t'} $ and $ \Delta A_{\mathbf{x},t} = \Delta x \Delta y :$
\begin{equation}\label{nrelvariation}
\frac{\sigma_n \left(r \right)|_{max}}{ \langle n \left( \mathbf{x_j}, \Delta t \right) \rangle} \sim \frac{4 \pi e^{1/4} D_o \Delta t}{ N_{\mathbf{x_j}, \Delta t} \Delta x \Delta y }
\end{equation}

Thus, the following procedure can be used to determine $ \Delta A_{\mathbf{x'},t'} ,$ $ \Delta t ,$ and $ N_{\mathbf{x'},t'} :$\\

\noindent a) Based on the nature of the transport problem under study, identify the spatial resolution required; this choice determines $ \Delta x ,$ $ \Delta y ,$ and $ \Delta A_{\mathbf{x}, t} = \Delta x \Delta y .$ 

\vspace{0.3cm}

\noindent b) Use Eq. (\ref{areaest}) to determine $ \Delta t ; $ for problems where $ \Delta x \approx \Delta y ,$ $ \Delta t = \sqrt{\Delta x} / D_o \approx \sqrt{\Delta y} / D_o .$

\vspace{0.3cm}

\noindent c) Choose a desired maximum relative (spatial, early-time) variation in RW's captured on $ \Delta A_{\mathbf{x}, t} , $ $ \sigma_n \left( \Delta t \right) / \langle n \left( \mathbf{x_j}, \Delta t \right) \rangle ,$ and solve Eq. (\ref{nrelvariation}) for $ N_{\mathbf{x'},t'} .$

\vspace{0.3cm}

\noindent Note 1: Since all RW's are independent, one can effectively increase the number of RW's launched from a given solution point, $(\mathbf{x},s)$ -  which, for single-swarm simulations, we denote as $N_{\mathbf{x},s}$ -  by relaunching another swarm of $M_{\mathbf{x},s}$ RW's from the same point. The number used in the relaunch, $M_{\mathbf{x},s},$ can be larger or smaller than $N_{\mathbf{x},s};$ the corresponding numbers of RW's captured in each interrogation area, $ \Delta A' , $ $ n\left( \mathbf{x}',s' \right) $ and $ m \left( \mathbf{x}',s' \right) ,$ respectively, can then be pooled to improve GF estimation accuracy. 

\vspace{0.3cm}

\noindent Note 2: As long as dimensional lengths are scaled by the either the short-time or long-time diffusion length scales, $ \sqrt{ D_o \Delta t }$ or $ \sqrt{D_o T } ,$ where $ T $ is the chosen solution time interval, this procedure applies for GF estimation procedures cast in either dimensional or nondimensional form.

\vspace{0.8cm}

\paragraph*{\textbf{S1 Appendix B}}
\label{S1 Appendix B}
{\bf  Properties of the respawning algorithm.}
In order to focus on the computational challenges associated with absorbing boundaries, in this subsection we limit attention to pure diffusion problems in which no advection takes place. Thus, the backward-time velocity field, $ \mathbf{b} ,$ is set to zero. The proposed respawning algorithm has three properties that appear to be essential to estimating low-error GF's within bounded domains: 

\vspace{0.3cm}

\noindent a) At all times and all locations in the solution domain, and with an error on the order of $ O \left( \Delta {x'}^2 , \Delta s' \right) , $ the algorithm conserves random walker mass, which, in turn, enforces the adjoint equation, (\ref{adjointeqnnew}).

In order to show that the respawning algorithm enforces the adjoint equation, label the location where the random walker having maximum weight at any instant $ {s'}_k > s , $ as $ {\mathbf{x}'}_j .$ Define a local $ \left( x" , y" \right)-\mathrm{coordinate} $ system, having origin at $ {\mathbf{x}'}_j ,$ such that the $ x"-\mathrm{axis} $ is aligned with the local gradient vector, $ \nabla' \Lambda , $ in the random walker mass density field, $ \Lambda = \Lambda \left( \mathbf{x}' , s' \right) .$ The continuous field, $  \Lambda \left( \mathbf{x}' , s' \right) , $ can be constructed, for example, by interpolation on the discrete instantaneous distribution of weighted RW's. 

The process of removing a maximum-weight random walker, having weight $ W_{max} \left( {s'}_k \right) ,$ located at $ {\mathbf{x}'}_j , $ and replacing it with two RW's, each of weight $ W_{max} \left( {s'}_k \right) /2 , $ and both placed at $ {\mathbf{x}'}_j , $ can be represented by the following expression:
\begin{equation}\label{spawneqn}
\frac{1}{2} \left[ \Lambda \left( {\mathbf{x}"}_j , {s'}_k | \mathbf{x}, s \right) +\Lambda \left( {\mathbf{x}"}_j , {s'}_k | \mathbf{x}, s \right) \right] - \Lambda \left( {\mathbf{x}"}_j , {s'}_k | \mathbf{x}, s \right) = 0
\end{equation}
random walker mass is clearly conserved since the left hand side of (\ref{spawneqn}) is identically zero. Taylor expanding the first term in the positive $ x"- $ direction, the second term in the negative $ x"- $ direction, and rearranging, we obtain:
\begin{equation}\label{spawneqntaylor}
\Lambda \left( {\mathbf{x}"}_j , {s'}_k | \mathbf{x}, s \right)= \frac{1}{2} \left[ \Lambda \left( {\mathbf{x}"}_{j+1} , {s'}_{k-1} | \mathbf{x}, s \right) +\Lambda \left( {\mathbf{x}"}_{j-1} , {s'}_{k-1} | \mathbf{x}, s \right) \right]  + O \left( {\Delta x"}^2 , \Delta s'  = \Delta {x"}^2/D \right)
\end{equation}
where, in the limit $ \Delta s' \rightarrow 0 ,$ $ D \Delta s' = \Delta {x"}^2 ,$ with $ \Delta {x"} = {x"}_{j+1}- {x"}_j = {x"}_j - {x"}_{j-1} ,$ $ \Delta s' = {s'}_k - {s'}_{k-1} , $ and where the order of discretization error is as shown.

Finally, expressing the left hand side of (\ref{spawneqntaylor}) as a temporal Taylor expansion about $ s' = {s'}_{k-1} ,$ and the two terms on the right as spatial expansions about $ {x"}_j ,$ we obtain the local mass diffusion equation, evaluated at $ \left( {\mathbf{x}'}_j , {s'}_{k-1} \right) , $ and expressed in the backward time direction
\begin{equation}\label{rwmassconservationeqn}
- \frac{\partial \Lambda}{\partial s'} = D \frac{\partial^2 \Lambda}{\partial {x"}^2}  + O \left( {\Delta x"}^2 , \Delta s' = \Delta {x"}^2/D \right)
\end{equation}
Note, that replacing the time index $ k $ with $ k+1 $ shifts evaluation of (\ref{rwmassconservationeqn}) to $ \left( {\mathbf{x}'}_j , {s'}_{k} \right) . $  

\vspace{0.3cm}

\noindent b) If we use the RW's to model physical Brownian particles, and the respawning algorithm to model some process by which the particles can undergo binary fission, we see that it enforces the diffusion equation (\ref{rwmassconservationeqn}) in the direction of the local random walker mass flux, $ \nabla \Lambda / \left| \nabla \Lambda \right| .$ Since physical particles, unlike numerical RW's, cannot both occupy the same point that their parent occupied, we would introduce a spatial offset for the respawned RW's. Thus, the fission of any given set of maximum-weight RW's, $ \{ \chi_{max,1} \left( {s'}_k \right) , \chi_{max,2} \left( {s'}_k \right) ,..., \chi_{max,M} \left( {s'}_k \right) \}, $ produces a mass gradient at each of the $ M  = M \left( s' \right) $ respawning locations, $ \chi_q \left( s' \right) , $ $ q=1,2,3,.., M \left(s' \right) , $ and the algorithm, in effect, smooths - on short (single-time-step) time-scales, $ \Delta {x'}^2 /D = \Delta s' $ - localized, and randomly-distributed (in space) removal of random walker mass.  Hence, the algorithm appears to be insensitive to the method used to reorder, at the end of every $ n $ time-steps, the weight array, $ \{W\} . $  Although this remains an open question, our results strongly support this view. 

\vspace{0.3cm}

\noindent c) As noted in the article, the third property associated with the respawning algorithm centers on the fact that it conserves the number of  RW's launched from the solution point, $ N_{\mathbf{x},s} . $ As shown in the article, this property ensures that accurate GF estimates can be obtained throughout any given time interval, $ 0 <  t' \leq T .$

\vspace{0.8cm}

\paragraph*{\textbf{S1 Appendix C}}
\label{S1 Appendix C}
{\bf Illustration of trial and error derivation of the adjoint equation and magic rule.}
A technical hurdle that likely limits widespread use of GF's, particularly among non-specialists, revolves around the origins of the adjoint equation governing the Green's function, and the so-called 'magic rule' \cite{barton}, the integral solution for the transport variable, $ \eta \left( \mathbf{x}, t \right) .$ 
Since textbooks and research articles appear to be universally mute on how these essential building blocks are obtained -  specifically, not showing how they can be simultaneously derived via straightforward trial and error - we illustrate trial and error determination of these the adjoint equation and magic rule, focusing for simplicity on fixed diffusivity  advection-diffusion problems, evolving in incompressible flow fields:
\begin{equation}\label{mainequationsuppinfo}
\frac{\partial\eta}{\partial t} + \mathbf{v}\cdot\nabla\eta - D \boldsymbol{\nabla^2} \eta= f(\mathbf{x}, t)
\end{equation}

The procedure incorporates four steps, the first of which is generally not highlighted: a) guess the form of the adjoint equation, b) multiply the candidate adjoint equation by  $ \eta \left( \mathbf{x}, t \right) , $ and the transport equation, Eq. (\ref{mainequationsuppinfo}) by the Green's function, $ G \left( \mathbf{x}, t | \mathbf{x'}, t') \right) ,$ c) integrate both equations in b) over the space-time hypersurface, $ \delta \Omega \times (0, t] ,$ formed by the product of the problem's spatial boundary, $ \delta \Omega , $ and the time axis, and d) taking the difference of the integral equations in c), use Green's theorems and trial and error alteration of the guessed adjoint equation to arrive at an integral solution for $ \eta \left( \mathbf{x}, t \right) ,$ stated strictly in terms of space-time surface integrals over $ \delta \Omega \times (0, t] ,$ the former involving $ G \left( \mathbf{x}, t | \mathbf{x'}, t') \right) , $ derivatives of $ G \left( \mathbf{x}, t | \mathbf{x'}, t') \right) , $ and known boundary and initial conditions.

In step a), it's generally good strategy to guess that the adjoint equation has the same generic form as the transport equation, Eq. (\ref{mainequationsuppinfo}), but with one or more signs reversed, and with the source term, $ f(\mathbf{x}, t),$ replaced by $ \pm \delta \left( \mathbf{x} -\mathbf{x'} \right) \delta \left( t - t' \right):$
\begin{equation}\label{guessadjointsuppinfo}
\frac{\partial G \left( \mathbf{x}, t | \mathbf{x'}, t') \right)}{\partial t'} + \mathbf{v}\cdot\nabla' G \left( \mathbf{x}, t | \mathbf{x'}, t') \right) + D \boldsymbol{{\nabla'}^2} G \left( \mathbf{x}, t | \mathbf{x'}, t') \right)=  \pm \delta \left( \mathbf{x} -\mathbf{x'} \right) \delta \left( t - t' \right)
\end{equation}
where: i) we have arbitrarily guessed that all signs on the left of the adjoint equation are positive, and ii) have fixed the solution point, $ \left( \mathbf{x} , t \right) .$ Choice ii): a) leads, in straightforward fashion, to an isolated term, $ \eta \left( \mathbf{x}, t \right) ,$ at the end of the procedure, and b) means that the  space-time integrals formed in step c) are taken over $ \mathbf{x'} $ and $ t' . $

Carrying out steps b) and c), we obtain:

\begin{align}
\label{magicsuppinfo1}
\int_{0}^{t+\epsilon} \int_{\Omega'} \eta \left( \mathbf{x'} , t' \right) \left[\frac{\partial G \left( \mathbf{x}, t | \mathbf{x'}, t' \right)}{\partial t'} + \mathbf{v}\cdot\nabla' G \left( \mathbf{x}, t | \mathbf{x'}, t' \right) + D \boldsymbol{{\nabla'}^2} G \left( \mathbf{x}, t | \mathbf{x'}, t' \right) \right] \hspace{.1cm} d \mathbf{x'} dt' - \nonumber \\
- \int_{0}^{t+\epsilon} \int_{\Omega'} G \left( \mathbf{x}, t | \mathbf{x'}, t' \right) \left[ \frac{\partial\eta \left( \mathbf{x'}, t' \right)}{\partial t} + \mathbf{v}\cdot\nabla' \eta \left( \mathbf{x'}, t' \right) - D \boldsymbol{\nabla'}^2 \eta \left( \mathbf{x'}, t' \right) \right]  d \mathbf{x'} dt' = \nonumber \\
\pm \int_{0}^{t+\epsilon} \int_{\Omega'}  \eta \left( \mathbf{x'} , t' \right) \mathbf{\delta \left( \mathbf{x} -\mathbf{x'} \right) \delta \left( t - t' \right)}  \hspace{0.1cm} d \mathbf{x'} dt' - \int_{0}^{t+\epsilon} \int_{\Omega'} G \left( \mathbf{x}, t | \mathbf{x'}, t' \right) f \left( \mathbf{x'}, t' \right) \hspace{0.1cm} d \mathbf{x'} dt'
\end{align}
where $ \mathbf{v} = \mathbf{v} \left( \mathbf{x'} , t' \right) ,$ and where the upper limit on the time integrals, $ t + \epsilon , $ is set just beyond $ t ; $ as shown in Eq. (\ref{magictimeterm}) below, following time integration, this trick suppresses the a priori unknown integral at $ t' = t+ \epsilon. $ In addition, it allows exposure of $ \eta \left( \mathbf{x}, t \right) $ from the first right-hand term in Eq. (\ref{magicsuppinfo1}).

Moving to step d) and focusing on Eq. (\ref{magicsuppinfo1}), use Liebniz rule on the fourth term on the left and combine with the first:
\begin{equation}
\label{magictimeterm}
\int_{\Omega'} \left[G \left( \mathbf{x}, t | \mathbf{x'}, t' \right) \eta \left( \mathbf{x'}, t' \right)   \right]_{0}^{t'+\epsilon} \hspace{0.1cm} d \mathbf{x'} - 2 \int_{0}^{t+\epsilon} \int_{\Omega'} \eta \left( \mathbf{x'} , t' \right) \frac{\partial G \left( \mathbf{x}, t | \mathbf{x'}, t' \right)}{\partial t'} \hspace{0.1cm} d \mathbf{x'} dt'
\end{equation}
Using the property that $ G \left( \mathbf{x}, t | \mathbf{x'}, t' \right) $ is undefined at for $ t' > t ,$ the first term simplifies to $ - \int_{\Omega'} G \left( \mathbf{x}, t | \mathbf{x'}, t'=0 \right) \eta \left( \mathbf{x'}, t'=0 \right)  \hspace{0.1cm} d \mathbf{x'} , $ where $ \eta \left( \mathbf{x'}, t'=0 \right) $ is the specified initial condition. Since the line (d=2) or area (d=3) integral in the second term of Eq. (\ref{magictimeterm}) must be suppressed, we see that the guessed plus sign on the first term in Eq. (\ref{guessadjointsuppinfo}) must be changed to a minus. 

Using similar steps on the two advective terms on the left side of Eq. (\ref{magicsuppinfo1} (terms 2 and 5) yields:
\begin{equation}
\label{magicadvecterm}
\int_0^{t+\epsilon} \int_{\Omega}' {\nabla'} \cdot \left( G \left( \mathbf{x}, t | \mathbf{x'}, t' \right) \eta \left( \mathbf{x'}, t' \right) \mathbf{v} \right) \hspace{0.1cm} d \mathbf{x'} dt' - 2 \int_0^{t+\epsilon} \int_{\Omega'}  G \left( \mathbf{x}, t | \mathbf{x'}, t' \right) {\nabla'} \cdot \left( \eta \left( \mathbf{x'}, t' \right) \mathbf{v} \right) \hspace{0.1cm} d \mathbf{x'} dt'
\end{equation}
where, for incompressible flow, $ {\nabla}' \mathbf{\cdot} \mathbf{v} = 0 ,$ allowing the $ \mathbf{v} $ in term 5, Eq. (\ref{magicsuppinfo1}) to be brought in front of the divergence operator: $  \mathbf{v}\cdot\nabla \eta \left( \mathbf{x'}, t' \right) = {\nabla}' \cdot \left( \mathbf{v} \eta \right) . $ Both area (d=2) or volume (d=3) integrals in Eq. (\ref{magicadvecterm}) are a priori unknown.  The first can be restated as a known/calculable line (d=2) or area (d=3) integral via the divergence theorem. The second, however, must again be suppressed by changing the guessed plus sign on term 2, Eq. (\ref{guessadjointsuppinfo}), to a minus.  Thus, the two terms in Eq. (\ref{magicadvecterm}) simplify to:
\begin{equation}
\label{magicadvectermnew}
\int_0^{t+\epsilon} \oint_{\delta \Omega_{advec}}  \left[ G \left( \mathbf{x}, t | \mathbf{x'}, t' \right) \eta \left( \mathbf{x'}, t' \right) \mathbf{v} \right] \mathbf{\cdot} \mathbf{\hat{n}'} \left( \mathbf{x'}, t' \right) \hspace{0.1cm} d \mathbf{x'} dt'
\end{equation}
where $ \delta \Omega_{advec}$ denotes the advective (inflow plus outflow) boundaries of $ \Omega ,$ and where $ \eta \left( \mathbf{x'}, t' \right) $ and $ \mathbf{v} \left( \mathbf{x'}, t' \right) $ are assumed known, at all times and locations on $ \delta \Omega_{advec}.$ In addition, $ \mathbf{\hat{n}'} \left( \mathbf{x'}, t' \right) $ is the  outward unit normal to $ \delta \Omega_{advec} .$ 
Finally, focusing on terms 3 and 6 in Eq. (\ref{magicsuppinfo1}), we are again lead to an a priori unknown surface/volume (d=2/3) integral:
\begin{equation}
\label{magicdiffnextra}
2 \int_0^{t+\epsilon} \int_{\Omega'}  {\nabla'} G \left( \mathbf{x}, t | \mathbf{x'}, t' \right) \mathbf{\cdot} {\nabla'} \eta \left( \mathbf{x'}, t' \right) \hspace{0.1cm} d \mathbf{x'} dt'
\end{equation}
In order to suppress this term, we change the plus sign on the third term in Eq. (\ref{guessadjointsuppinfo}) to a negative.  Using the divergence theorem and integration by parts, introducing the three sign changes above into the guessed adjoint equation, Eq. (\ref{guessadjointsuppinfo}), labeling the functions associated with specified Dirichlet, Neumann, and initial conditions on $ \eta \left( \mathbf{x}, t \right) ,$ respectively, as $ g \left( \mathbf{x} , t \right) ,$ $ h \left( \mathbf{x} , t \right) ,$ and $ \phi \left( \mathbf{x} \right) ,$ and 
choosing the plus sign on the first term of the right side of Eq. (\ref{magicsuppinfo1}), we finally obtain the magic rule: 
\begin{align}
\label{magicsuppinfo2}
\eta \left( \mathbf{x}, t \right) =  \int_{0}^{t} \int_{\Omega'} G \left( \mathbf{x}, t | \mathbf{x'}, t' \right) f \left( \mathbf{x'}, t' \right) \hspace{0.1cm} d \mathbf{x'} dt' + \int_{\Omega'} G \left( \mathbf{x}, t | \mathbf{x'}, t'=0 \right) \eta \left( \mathbf{x'}, t'=0 \right)  \hspace{0.1cm} d \mathbf{x'} - \hspace{0.2cm} \nonumber \\
- D \int_{0}^{t} \int_{\Omega_D'} h \left( \mathbf{x'} , t' \right) \boldsymbol{{\nabla'}} G \left( \mathbf{x}, t | \mathbf{x'}, t' \right) \cdot \mathbf{n'}  \hspace{.1cm} d \mathbf{x'} dt' + D \int_{0}^{t} \int_{\Omega_N'} G \left( \mathbf{x}, t | \mathbf{x'}, t' \right)  \boldsymbol{{\nabla'}} \eta \left( \mathbf{x'} , t' \right)  \cdot \mathbf{n'}  \hspace{.1cm} d \mathbf{x'} dt' - \hspace{0.1cm} \nonumber \\
- \int_{0}^{t} \int_{\Omega_F'} G \left( \mathbf{x}, t | \mathbf{x'}, t' \right) \eta \left( \mathbf{x'} , t' \right) \mathbf{v} \left( \mathbf{x'} , t' \right) \cdot \mathbf{n'} \hspace{.1cm} d \mathbf{x'} dt' \hspace{.6cm}
\end{align}
where $ \Omega'_D ,$ $ \Omega'_N ,$ and $ \Omega'_F $ correspond, respectively, to the portions of the domain boundary, $ \delta \Omega' ,$ on which Dirichlet, Neumann, and inflow/outflow boundary conditions are imposed, and where $ h \left( \mathbf{x'}, t' \right) $ is the spatially and temporally specified value of $ \eta \left( \mathbf{x'}, t' \right) $ on the Dirichlet boundary.  

\vspace{0.8cm}

\paragraph*{\textbf{S1 Appendix D}}
\label{S1 Appendix D}
{\bf Extraction of the transport equation, Eq. (\ref{mainequationfwd}), associated with the Green's function in test 3.}  
The transition density estimated in test 3, governed by the Fokker-Planck equation, Eq. (\ref{fptest3}), also satisfies the backward Kolmogorov equation:
\begin{equation}\label{bketest3}
\mathbf{D'} \mathbf{:} {\nabla'}^T \otimes \nabla' p + \mathbf{b} \cdot \nabla' p +  \frac{\partial p}{\partial t'} = 0,
\end{equation}
where $ \mathbf{D'} = \mathbf{D'} \left( \mathbf{x'} , t' \right) ,$ $ \nabla' ,$ and $ \mathbf{b} = \mathbf{b} \left( \mathbf{x'} , t' \right) = - \mathbf{v} \left( \mathbf{x'} , t' \right) ,$ are expressed in terms of variable 'backward' variables, $ \mathbf{x'} $ and $ t' .$ 

The corresponding adjoint equation is given by: 
\begin{equation}\label{propeqnapp}
\mathbf{D'} \mathbf{:} {\nabla'}^T \otimes \nabla' K + \mathbf{b} \cdot \nabla' K +  \frac{\partial K}{\partial t'} = 0,
\end{equation}
where, again, the Green's function is given by:
\begin{equation}\label{propandgf2}
G \left( \mathbf{x} , t | \mathbf{x'},t' \right) = H \left( t -t' \right) K \left( \mathbf{x} , t | \mathbf{x'},t' \right)
\end{equation}
Given Eq. (\ref{propeqnapp}), the associated transport equation is given by:
\begin{equation}\label{mainequationapp}
\mathbf{D} \mathbf{:} {\nabla}^T \otimes \nabla  \eta - \mathbf{b}\cdot\nabla\eta - \frac{\partial\eta}{\partial t} = -f(\mathbf{x}, t),
\end{equation}
where, in Eq. (\ref{mainequationfwd}), we write $ \mathbf{b} $ as $ \mathbf{v_a} .$


\bibliographystyle{elsarticle-num}






\end{document}